\input amssym.def
\input amssym.tex
\magnification 1200
\hoffset= 2pc
\hsize=30pc
\vsize=45pc
\font\sectionfont=cmbx12 at 12pt

\def\NSF98{\footnote*{This research was supported in part by NSF grant
 DMS-9803593.}}
\def\widetilde{\mathaccent"0365 }

\def\align{\eqalign}
\def\fk{\frak k}
\def\ft{\frak t}
\def\fg{\frak g}
\def\tcalC{{\tilde\calC}}
\def\ot{\otimes}
\def\pr{\prime}
\def\ze{\zeta}
\def\a{\alpha}
\def\b{\beta} 

\def\g{\gamma}
\def\O{\Omega}
\def\o{\omega}
\def\un{\underline}
\def\und{\underbar}
\def\th{\theta}
\def\la{\lambda}
\def\La{\Lambda}
\def\inf{\infty}
\def\s{\sigma}
\def\na{\nabla}
\def\noblackbox{\overfullrule=0pt}
\def\jlb{J-L. Brylinski}
\def\Sp{Springer Verlag}
\def\bu{\bullet}
\def\BC{{ \Bbb C}}
\def\BQ{{ \Bbb Q}}
\def\BR{{ \Bbb R}}
\def\BT{{ \Bbb T}}
\def\BZ{{ \Bbb Z}}
\def\calC{{\cal C}}
\def\calD{{\cal D}}

\def\calL{{\cal L}}
\def\calU{{\cal U}}
\def\calV{{\cal V}}
\def\calO{{\cal O}}
\def\lan{\langle}
\def\ran{\rangle}
\def\pha{\phantom}
\def\display{\displaystyle}

\def\longto{\longrightarrow}
\def\mapright#1{\smash
{\mathop{
\longto}\limits^{#1}}}
\def\qed{\vrule height4pt 
width3pt depth2pt}
\def\mapup#1{\Big\uparrow\rlap{$\vcenter
{\hbox{$\scriptstyle#1$}}$}}
\def\mapdown#1{\Big\downarrow
\rlap{$\vcenter
{\hbox{$\scriptstyle#1$}}$}}

\def\theo#1#2{\vskip 1pc
  {\bf   Theorem\ #1.} {\it #2}
\vskip 1pc}
\def\lem#1#2{\vskip 1pc{\bf 
\ Lemma #1.}
  {\it #2}\vskip 1pc}
\def\cor#1#2{\vskip 1pc{\bf 
\ Corollary\ #1.}
  {\it #2}\vskip 1pc}
\def\prop#1#2{\vskip 1pc{\bf 
 Proposition \ #1.} 
{\it #2}}

\def\proof#1{\vskip 1pc {\bf  
Proof.} {\rm #1} 
\vskip 1pc}

\def\sec#1{\vskip 1.5pc\noindent
{\hbox 
{{\sectionfont #1}}}\vskip 1pc}

\def\min{\setminus}
\noblackbox

\centerline{\bf GERBES ON COMPLEX REDUCTIVE LIE GROUPS}
\vskip 1.2pc

\centerline{Jean-Luc Brylinski \NSF98}
\vskip 2pc

\bf Abstract \rm  
Let $K$ be a compact Lie group with complexification $G$.
We construct by geometric methods a conjugation invariant gerbe on
$G$;  this then gives by restriction an invariant gerbe on $K$.
Our construction works for any choice of level. When 
$K$ is simple and simply-connected, the level is just an integer as
usual. For general $K$, the level is a  bilinear form $b$ on a Cartan
subalgebra where $b$ satisfies a quantization condition. 

The idea of our construction is to  first construct a gerbe on 
the Grothendieck manifold of pairs $(g,B)$ where $g\in G$ and $B$ is a
Borel subgroup containing $g$. Then the main work is to descend that
gerbe to $G$. There is an interesting torsion phenomenon in that the  
restriction of the gerbe to a semisimple orbit is not always trivial. 

The paper starts with a discussion of gerbe data and of 
gerbes as geometric objects (sheaves of
groupoids); the relation between the two approaches is presented. The  
Appendix on equivariant gerbes  discusses   both points of view. 
\vskip 1pc \rm 

\vskip .13 pc

\sec{0. Introduction}
Gerbes are higher analogs of  bundles. The gerbes we consider here
are the so-called DD gerbes (for Dixmier and Douady), and they are analogs
of line bundles.
General gerbes were introduced by Giraud [{\bf Gi}] for the purpose
of his theory of degree $2$ non-commutative sheaf cohomology.
The differential geometry of DD gerbes was developed in the book [{\bf Br1}],
where the analogs of connection and curvature were introduced.
In particular, the curvature $\O$ of a gerbe is a $3$-form, called the
$3$-curvature. This $3$-form is quantized in the same way as the curvature
of a line bundle is quantized. Gerbes are
classified by their
$3$-curvature up to so-called flat gerbes.
The methods of [{\bf Br1}] are sheaf-theoretic.
Gerbes can be obtained from an open covering $(V_\a)$ and from
line bundles $(\La_{\a\b})$ on overlaps $V_{\a\b}$ together with 
trivializations of  $\La_{\a\b}\ot\La_{\b\g}\ot \La_{\g\a}$
over $V_{\a\b\g}$. This  point of view was first introduced in [{\bf Br1}]
for the example of $S^3$ covered by two balls, in connection with the magnetic
monopole. It was developed systematically by Chatterjee [{\bf Ch}]
and Hitchin [{\bf Hi}]. The relation with (smooth and holomorphic) Deligne
cohomology is discussed in [{\bf Br1}].

One of the first instances of gerbes occurs on a simple simply-connected compact Lie group
$K$. We have the well-known normalized Chern-Simons $3$-form $\nu$ on $K$,
and it is proved in [{\bf Br1}] that there is a completely \und{canonical}
gerbe $\calC$ on $K$ with $3$-curvature $\O=2\pi i\nu$. We gave an explicit
construction of this gerbe in [{\bf Br1}] using the path-fibration $PK\to K$ and the
central extension $\widetilde{\O K}$ of the based loop group $\O K$. 
One can view this gerbe as a geometric realization of a class in smooth Deligne
cohomology (or equivalently, the group of differential characters of Cheeger-Simons).

However for a number of reasons it is very interesting to have a geometric
construction of the gerbe $\calC$ on $K$ which only invokes finite-dimensional
geometry. One reason is that one fundamental use of gerbes is to construct line
bundles over the free loop space and in particular central extensions
of loop groups [{\bf Br1}] [{\bf Br-ML}]. This is a good reason to want a
description of the canonical gerbe
$\calC$ which does not use the central extension of $\O K$.
Another motivation is the  theory of group valued moment maps
introduced by Alekseev, Malkin and Meinrenken [{\bf A-M-M}],
which is closely related to the canonical gerbe on a compact Lie group.
This theory really belongs to finite-dimensional geometry. Indeed
Alekseev, Meinrenken and Woodward show how to use group valued moment maps to study
hamiltonian actions of loop groups with proper moment map and obtain index formulas [{\bf
A-M-W}]. There is a description of a gerbe over the complexification $G$ of $K$ which uses
the Bruhat decomposition of $G$ into double cosets $BwB$ where $B$
is a Borel subgroup of $G$. Such a description is implicit in the paper
[{\bf B-D}] where certain cohomology classes are constructed in algebraic K-theory. 
According to [{\bf Br2}], these classes yield holomorphic gerbes. 
The ensuing gerbe on $G$ is manifestly equivariant under left translation by $B$. The
point of view here is that one gets stronger results by constructing a holomorphic gerbe
on $G$; given a gerbe over $G$ one can simply restrict it to $K$ to get a gerbe on $K$.

Here we also use the complexification $G$ and construct a holomorphic gerbe
over $G$. Our method is to use an auxiliary manifold $\tilde G$, called the
Grothendieck manifold, which is the set of pairs $(g,B)$ where  $B\subset G$ is a Borel
subgroup and
$g\in B$. We have a projection map $q:\tilde G\to G$ whose fiber over $g$ is the variety
of Borel subgroups which contain $g$. Our method is as follows: first  we start from some
combinatorial data, which is an element $b\in X^*(T)\ot X^*(T)$, where $T$ is a maximal
torus in $G$ and $X^*(T)$ is the character group of $T$. From this combinatorial data
we easily construct some gerbe
$\tcalC$ over $\tilde G$. The idea is that there are natural mappings from $\tilde G$ to 
the complete flag manifold $G/B$ and to $T$. Thus one can pull-back to $\tilde G$
characters of $T$ and equivariant line bundles over $G/B$. Then there is a natural
cup-product construction which creates a gerbe $\tcalC$ over $\tilde G$, which has an explicit
description by gerbe data.
Then we want to descend $\tcalC$ to a gerbe on
$G$. The idea is the following. First over the open set $G^{reg}$ of regular semisimple
elements, the mapping $\tilde G \to G$ restricts to give a Galois covering whose
Galois group is the Weyl group $W$. If we assume that $b$ is $W$-invariant, then we
can construct a gerbe over
$G^{reg}$ by descent. Now we need to extend this gerbe to all of $G$.
This is done in several steps, assuming that $b(\check\a,\check\a)$ is even for all
coroots $\a$; first we extend $\calC$ across some divisors, by reducing the problem to the
case of
$SL(2,\BC)$ (\S 3 and 4). Then (\S 5) we have to handle codimension $2$ subvarieties by
cohomological methods. We also obtain a
$0$-connection on our gerbe $\calC$. We don't have explicit data for the gerbe $\calC$,
because we have to use abstract processes to descend it and extend it. Instead
the gerbe is constructed as a sheaf of groupoids. The construction
of the next differentiable structure (a so-called
$1$-connection) is still conjectural.

Our construction is in fact more general in that we consider an arbitrary
complex reductive group $G$ (for instance $GL(n,\BC)$), and we construct
holomorphic gerbes on
$G$ which are equivariant under some auxiliary group $H$, which is only equal to
$G$ up to center. The combinatorial data then becomes a tensor in
$X^*(S)\ot X^*(T)$, where $S$ is a maximal torus in $H$.

It is interesting to restrict the equivariant gerbe $\calC$ to the conjugation orbit
of some $g\in G$. The obstruction to the triviality of the restricted gerbe
is a central extension of the centralizer group. This
is discussed in \S 7 for semisimple elements, where an explicit cocycle
is given for the extension. Outside of
$SL(n,\BC_)$ the   central extensions can be non-trivial, and this is an obstacle to
constructing equivariant gerbe data for our gerbe in general, as Alekseev, Meinrenken and
Woodward were able to do for $SL(n)$.

We have included in the first section an exposition of the DD gerbes
and their differential geometry. First we expose the basics of the gerbe data in the
sense of Chatterjee [{\bf Ch}] and Hitchin [{\bf Hi}]. This has the advantage
of being very concrete and explicit. However for this paper we need to use the more
general and abstract formalism of [{\bf Br1}] based on sheaves of groupoids so we
discuss that next, including the differential geometry, and we give a brief
comparison between the two approaches. We intend to soon write a more detailed
account of this in an expository paper.

The equivariant gerbes are discussed in an Appendix, starting from the gerbe data
approach and then leading to the sheaves of groupoids.

The second section gives a discussion of the geometry of the manifold $\tilde G$
and of its very interesting cohomology. Particularly noteworthy is the fact
that we have $H^*(G)=H^*(\tilde G)^W$ for cohomology with rational coefficients.

It is a pleasure to thank Anton Alekseev, Pierre Deligne, Johannes Huebschmann, Eckhard
Meinrenken and
 Victor Guillemin for many useful discussions. In particular Alekseev and Meinrenken
informed me about their recent works on hamiltonian reduction for loop group actions
and told me about their construction of a gerbe on the unitary group.

\sec{1. Gerbe data versus gerbes}

First we discuss gerbe data in the sense of Chatterjee [{\bf Ch}] and Hitchin [{\bf
Hi}]. Gerbe data on  a manifold $M$ consists of the following:

1) an open covering $(V_\a)$ of $\tilde M$;

2) a family of line bundles $\La_{\a\b}\to V_{\a\b}$ such that
$\La_{\a\a}=\bf 1$ is the trivial line bundle

3) an isomorphism $u_{\a\b}:\La_{\a\b}^{\ot-1}\tilde{\to}\La_{\b\a}$
such that, viewed as a section of $\La_{\a\b}\ot\La_{\b\a}$, $\phi_{\a\b}$ is
symmetric in $(\a,\b)$

4) for each $\a,\b,\g$,  a non-vanishing section $\th_{\a\b\g}$
of the tensor product $\La_{\a\b}\ot\La_{\b\g}\ot
\La_{\g\a}$ over $V_{\a\b\g}$ satisfying the cocycle condition
$$\th_{\b\g\delta}\ot\th_{\a\g\delta}\ot\th_{\a\b\delta}
\ot\th_{\a\b\g}=1$$ over $V_{\a\b\g\delta}$. Note this makes sense as the
quadruple tensor product is easily seen to be a section of a trivial line
bundle. We require that $u_{\a\b}\ot u_{\b\g}\ot u_{\g\a}$ transforms
$\th_{\a\b\g}$ into $\th_{\g\b\a}^{\ot-1}$.

It is of course important to know when two gerbe data are equivalent.
The notion of equivalence is generated by two types of operations.
The first operation is simply restriction to a finer open covering.
The second type of operation is a kind of gauge transformation. There
are actually two kinds of gauge transformations acting on gerbe data
for a fixed open covering $(V_\a)$. Firstly, while keeping the line bundles
$\La_{\a\b}$ fixed, we can pick smooth functions $h_{\a\b}:V_{\a\b}\to\BC^*$
and change $\th_{\a\b\g}$ to $h_{\b\g}h_{\a\g}^{-1}h_{\a\b}\th_{\a\b\g}$.
Secondly, we can introduce auxiliary line bundles $E_{\a}$ over $V_\a$
and change the line bundle $\La_{\a\b}$ to $\La'_{\a\b}=E_\a\ot\La_{\a\b}\ot
E_\b^{\ot-1}$; then the triple tensor product 
$\La'_{\a\b}\ot\La'_{\b\g}\ot
\La'_{\g\a}$ is canonically isomorphic to
$\La_{\a\b}\ot\La_{\b\g}\ot
\La_{\g\a}$, so we take the same $\th_{\a\b\g}$, now viewed as a section
of $\La'_{\a\b}\ot\La'_{\b\g}\ot\La'_{\g\a}$.
Then two gerbe data are equivalent if they are related by a sequence of the
various operations we just discussed.

It is useful to single out the criterion for a gerbe data $(\La_{\a\b},\th_{\a\b\g})$
to be trivial. This means that we can find line bundles $E_\a$ over $V_\a$,
and isomorphisms $E_\a^{\ot-1}\ot E_\b \tilde{\to}\La_{\a\b}$
over $V_{\a\b}$ such that $u_{\a\b}$ is the obvious isomorphism and
$\th_{\a\b\g}$
corresponds to the section
$1$ of the trivial line bundle
$$[ E_\b^{\ot-1}\ot E_\g]\ot [ E_\a^{\ot-1}\ot E_\g]
^{\ot-1}\ot [E_\a^{\ot-1}\ot E_\b].$$

Very often,  as  in [{\bf Hi}], we will not include $u_{\a\b}$ explicitly in the gerbe data.

However perhaps it is clearer to view the notion of equivalence of gerbe
data in terms of tensor products of gerbe data, which we discuss next.

Given two gerbe data on $M$, we can first of all refine the open coverings
used to define them so that they become the same. Then we are in the situation
of an open covering $(V_\a)$ and two gerbe data
$(\La_{\a\b},\th_{\a\b\g})$ and $(M_{\a\b},\rho_{\a\b\g})$. Then their tensor
product is the data $(\La_{\a\b}\ot M_{\a\b},\th_{\a\b\g}\ot
\rho_{\a\b\g})$. The inverse of the gerbe data $(\La_{\a\b},\th_{\a\b\g})$
is simply $(\La_{\a\b}^{\ot-1},\th_{\a\b\g}^{\ot-1})$.

Then two gerbe data $(\La_{\a\b},\th_{\a\b\g})$ and $(M_{\a\b},\rho_{\a\b\g})$
over the same open covering are equivalent if the data
$(\La_{\a\b}\ot
M_{\a\b}^{\ot-1},\th_{\a\b\g}\ot\rho_{\a\b\g}^{\ot-1})$
is a trivial gerbe data.

One nice feature of gerbe data is that it is very easy to pull them back by an
arbitrary smooth mapping $f:Y\to M$: one simply introduces the open covering
$f^{-1}V_\a$ of $Y$, together with the pull-back line bundles
$f^*\La_{\a\b}$ and the pull-backs of the $\th_{\a\b\g}$.

Then there is the notion of a $0$-connection. This is given by a family
of connections $D_{\a\b}$ on the $\La_{\a\b}$, such that $u_{\a\b}$ is compatible
with the connections and 
$\th_{\a\b\g}$ is horizontal with respect to the tensor product
connection on $\La_{\a\b}\ot\La_{\b\g}\ot
\La_{\g\a}$. Our two kinds of gauge transformations
for gerbe data  can be extended to the $0$-connections. The first type
of extended gauge transformation depends on functions $h_{\a\b}:V_{\a\b}\to\BC^*$ as
before; besides transforming $\th_{\a\b\g}$ into
$h_{\b\g}h_{\a\g}^{-1}h_{\a\b}\th_{\a\b\g}$, the connection $D_{\a\b}$ is transformed
into $D_{\a\b}+d \log(h_{\a\b})$. The second type depends on auxiliary line bundles
$E_\a$ equipped with connections
$\na_\a$; then when $\La_{\a\b}$ is changed to $\La'_{\a\b}=E_\a\ot\La_{\a\b}\ot
E_\b^{\ot-1}$, the new line bundle acquires the tensor product
connection $\na_\a+D_{\a\b}-\na_\b$.

Next we have the notion of a $1$-connection. This consists of $2$-forms $F_\a$ on
$V_\a$ such that
$$F_\b-F_\a=Curv(D_{\a\b}),$$
where $Curv(D_{\a\b})$ denotes the curvature of $D_{\a\b}$. Then only the second kind of
gauge transformation acts on the $2$-forms, transforming $F_\a$ into $F_\a+Curv(D_\a)$.

Given a $1$-connection, we obtain a global closed $3$-form $\O$ such that
$\O_{/V_\a}=d F_\a$. This $3$-form is called the $3$-curvature.

\theo{1}{\rm [{\bf Br1}\rm] The $3$-curvature $\O$ is quantized, i.e., the periods
of $(2\pi i)^{-1}\O$ are integral. Conversely, every quantizable
$3$-form occurs as the $3$-curvature of some gerbe data.}

The cohomology class is (when $H^3(M,\BZ)$ has no torsion) the only obstruction
to the triviality of the gerbe data.

Let us consider the question of classifying the equivalence classes of gerbe data
over $M$. By picking a fine enough open covering, we may assume that all line bundles
$\La_{\a\b}$ are trivial, in which case we have $\th_{a\b\g}=g_{\a\b\g}$
where $g_{\a\b\g}:V_{a\b\g}\to\BC^*$ is a \v Cech $2$-cocycle with values in the
sheaf $\un{\BC}^*$ of smooth $\BC^*$-valued functions. Then we have

\prop{1}{The equivalence classes of gerbe data over $M$ are classified by the \v
Cech cohomology group $H^2(M,\un{\BC}^*)=H^3(M,\BZ)$.}

The point here is that when the $\La_{\a\b}$ are trivial, the only gauge transformation we
can still use on gerbe data is $g_{\a\b\g}\mapsto
g_{\a\b\g}h_{\b\g}h_{\a\g}^{-1}h_{\a\b}$, which amounts to  multiplying the \v
Cech cocycle $(g_{\a\b\g})$ by a coboundary.

We now consider gerbes with $0$-connections: as $\La_{\a\b}$ is the trivial line
bundle, the connection $D_{\a\b}$ is simply a $1$-form $A_{\a\b}$.
These $1$-forms satisfy
$$A_{\a\b}+A_{\b\g}+A_{\g\a}=d\log (g_{\a\b\g}).$$ 
This means that the pair $(g_{\a\b\g},A_{\a\b})$ is a \v Cech cocycle with
respect to the complex of shaves $\un{\BC}^*\mapright{d\log}\un A^1_M$,
where $\un A^1_M$ is the sheaf of $1$-forms on $M$. The gauge transformation
$$(g_{\a\b\g},A_{\a\b})\mapsto (g_{\a\b\g}h_{\b\g}h_{\a\g}^{-1}h_{\a\b},
A_{\a\b}+d\log h_{\a\b})$$
corresponds to a \v Cech coboundary. Thus gerbe data with $0$-connection are
parameterized up to equivalence by the \v Cech hypercohomology group
$H^2(M,\un{\BC}^*\mapright{d\log} \un A^1_M)$.
This is the smooth Deligne cohomology group $H^3(M,\BZ(2)^{\inf}_D)$; see
[{\bf Br1}] and [{\bf D-F}] for smooth Deligne cohomology groups
and their relation to classical field theory

Similarly given gerbe data with $0$-connection and $1$-connection, one obtains \v
Cech cocycle $(g_{\a\b\g},A_{\a\b},F_\a)$ which is unique up to a coboundary.
Thus gerbe data equipped with $0$-connection and $1$-connection
are classified by the \v Cech hypercohomology group
$H^2(M,\un{\BC}^*\mapright{d\log} \un A^1_M\mapright{d}\un A^2_M)$.
Again this is a smooth Deligne cohomology group.

As concrete and flexible as the formalism of gerbe data is, it does
have some limitations. First of all, it may hard to manipulate in situations
where it is not possible to work with a unique open covering. For instance,
when a group action is present and we look for equivariant gerbes, it will
not always be possible to work with invariant open sets (see Appendix).
Secondly, it leaves open the question of a geometric interpretation
of what lies behind the gerbe data. In other words, if we view gerbes as the higher
analog of line bundles, what takes the place of the total space of a gerbe?

Although there is nothing like a manifold intrinsically associated to a gerbe, there
are plenty of geometric objects which arise when we think of the gerbe data
$(\La_{\a\b},\phi_{\a\b},\th_{\a\b\g})$ as instructions to solve a geometric
construction problem. The problem is to trivialize the gerbe, which means
to find line bundles $E_\a$ over $V_\a$,
and isomorphisms $\psi_{\a\b}: E_\a^{\ot-1}\ot E_\b\tilde{\to}\La_{\a\b}$
over $V_{\a\b}$ such that $\th_{\a\b\g}$ corresponds to the section $1$
of the trivial line bundle 
$$[ E_\b^{\ot-1}\ot E_\g]\ot [ E_\a^{\ot-1}\ot E_\g]
^{\ot-1}\ot [E_\a^{\ot-1}\ot E_\b].$$ The problem is that if the gerbe data is not
trivial, there is no solution, so we are dealing with an insoluble construction
problem. However, we can solve the problem at least locally, over some small enough
open set $U$ of $M$. This means that the line bundles $E_\a$ are defined over
the intersection $U\cap V_\a$ and the isomorphisms $\phi_{\a\b}$ are defined
over $U\cap V_{\a\b}$.

Then what takes the place of a total space for a gerbe is the collection $\calC_U$
of all such data $(E_\a,\psi_{\a\b})$ as $U$ ranges over all open sets in $M$.
In this collection, it would be a fatal mistake to identify all isomorphic
objects. Indeed, it is essential to keep precise track of how objects
of $\calC_U$ restrict to smaller open sets and conversely how we can glue together
objects of $\calC_{U_i}$ into an object of $\calC_U$, when $(U_i)$ is an open
covering of $U$. For thus purpose, we need to take full account of the inner
structure of $\calC_U$: it is not just a collection of objects, but we have the
notion of isomorphism of objects. Thus each $\calC_U$ is a groupoid (i.e., a category
where each morphisms is invertible), and we can view the collection
of all $\calC_U$ as a sheaf of groupoids over $M$. This means that for $V\subset U$
we have a restriction functor $\calC_U\to \calC_V$ (usually denoted on objects of
$\calC_U$  by $P\mapsto P_{/V}$), and if $(U_\a)$ is an open covering of some open
set $U$, it amounts to the same to give an object of $\calC_U$ or to give objects
$P_\a$ of $\calC_{U_\a}$ together with transition isomorphisms
$$\psi_{\a\b}:[P_\b]_{/U_{\a\b}}\tilde{\to}[P_\a]_{/U_{\a\b}}$$
with $\psi_{\a\b}=\psi_{\b\a}^{-1}$, which satisfy the cocycle condition
 $\psi_{\a\b}\psi_{\b\g}\psi_{\g\a}=Id$.

We refer to the book [{\bf Br1}] for a full discussion of sheaves of groupoids,
which are essentially the stacks introduced by Grothendieck around 1961. A sheaf of
groupoids is called a \und{Dixmier-Douady gerbe} (or DD-gerbe)  if it satisfies the
following:

1) all objects of $\calC_U$ are locally isomorphic;

2) for any $x\in M$, there is an open set $U$ containing $x$ such that
$\calC_U$ is non-empty;

3) for any object $P$ of $\calC_U$, the automorphisms of $P$ are exactly
the smooth functions $U\to\BC^*$.

We have seen how gerbe data leads to such DD-gerbes. Hereafter, when we talk about gerbes,
we always have in mind DD gerbes.
Next we explain briefly how a gerbe leads to gerbe data. Given two objects
$Q,R$ of some $\calC_U$, the isomorphisms between $Q$ and $R$ are the sections
of some $\BC^*$-bundle $\un{Isom}(Q,R)$. This follows from axioms 1) and 3).
Using axiom 2) we can find an open covering
$(V_\a)$ and objects $P_\a\in\calC_{V_\a}$. Then over $V_{\a\b}$ we have the line
bundle
$\La_{\a\b}$ associated to the $\BC^*$-bundle $\un{Isom}(P_\b,P_\a)$.
Then the composition law for isomorphisms immediately yields the isomorphism
$u_{\a\b}$ as well as a non-vanishing section
$\th_{\a\b\g}$ of $\La_{\a\b}\ot\La_{\b\g}\ot
\La_{\g\a}$. The notion of gerbe gives a geometric explanation for the two types
of gauge transformations on gerbe data. The first type, corresponding to functions
$h_{\a\b}:V_{\a\b}\to\BC^*$ amounts to keeping the objects $P_\a$ and the line bundles
$\La_{\a\b}$ fixed but
changing the isomorphisms $\phi_{\a\b}:\La_{\a\b}\min {\rm ~(
0-section)~}\tilde{\to}\un{Isom}(P_\b,P_\a)$ by multiplication by $h_{\a\b}$.
The second type, given by line bundles $E_\a$, consists in changing the
objects $P_\a$ to $P_\a\ot E_\a$ by twisting  so that
$\La_{\a\b}=\un{Isom}(P_\b,P_\a)$ is changed to $\un{Isom}(E_\b\ot P_\b,E_\a\ot P_\a)=E_\a\ot
\La_{\a\b}\ot E_\b^{-1}$.

 See [{\bf Br1,\rm  Chap. 4 and 5}] for a detailed discussion of gerbes
and degree $2$ cohomology as well as many  examples.

We have seen that gerbes and gerbe data (up to gauge equivalence)
are essentially equivalent notions, and one can go back and forth between the two.
It is
interesting to note that geometric situations where gerbe occur usually lead directly to a
gerbe, and only secondarily to gerbe data, after some auxiliary choices are made. This is the
case for a principal
$G$-bundle $Q\to M$, when we are given a central extension $\tilde G$ of $G$ by $\BC^*$. Then
the objects of
$\calC_U$ are the principal $\tilde G$-bundles over $U$ which lift $Q_{/U}$. The gerbe data
appear when we pick liftings $\tilde Q_\a\to V_\a$; then we have $\tilde Q_\b=\tilde Q_\a\ot
\La_{\a\b}$ for a unique line bundle $\La_{\a\b}$, and the construction automatically yields
the section $\th_{\a\b\g}$. The case where $\tilde G$ is the group of invertible
linear transformations in a Hilbert space, and $G$ is the projective linear
group, is studied in depth in the article of Dixmier and Douady [{\bf D-D}].
This explains the name: Dixmier-Douady gerbes. The book [{\bf Br1}] presents
another example, where $\tilde G$ is a generalized Heisenberg group in the sense
of Weil.

The trivial gerbe is still interesting to think about, as it has lots of structure
and is the local model for non-trivial gerbes. The objects of $\calC_U$ are line
bundles over $U$, and isomorphisms are line bundle isomorphisms. Equivalently, we
can think of the objects of the trivial gerbe as being $\BC^*$-bundles.
As soon as a gerbe has an object over some open set, its restriction to that
open set becomes trivial. More precisely, we have:

\prop{2}{(1)Given a DD gerbe $\calC$ over $M$, for each open set $U$ of $M$,
and for any given object $P$ of $\calC_U$ we have an equivalence of categories
between $\calC_U$ and the category of $\BC^*$-bundles over $U$, given on objects by
$Q\mapsto \un{Isom}(P,Q)$.
\vskip .04 in
(2) If $\calC_U$ is not empty, its objects are classified up to isomorphism by the
group $H^2(U,\BZ)$.}

Thus for any line bundle $L\to U$ and for any object $P$ of $\calC$, we have a
well-defined twist $P\ot L$ of $P$ by $L$ such that $\un{Isom}(P,P\ot L)$ is the
$\BC^*$-bundle associated to $L$. This can be constructed by explicit gluing as follows.
Let $(g_{\a\b})$ be some transition functions for $L$, with respect to some open
covering $U_\a$ of $U$. Then $g_{\a\b}$ gives an automorphism of $P$ over $U_{\a\b}$,
which satisfies the cocycle condition, hence can be used to construct an object
of $\calC_U$.

The notion of equivalence of gerbes is also very natural from the point of view of
sheaves of groupoids.
Thus if $\calC$ and $\calD$ are two DD-gerbes on $M$, an equivalence
$\phi:\calC\to\calD$ of gerbes is a family of functors $\calC_U\to \calD_U$
(for $U$ an arbitrary open set in $M$)
which is compatible with restriction to smaller open sets. Given such $\phi$,
each $\phi_U$ is an equivalence of categories by Proposition 2 (2).

We next discuss the tensor product $\calC\ot\calD$ of two DD gerbes.
The tensor product  gerbe is such that if $P$ is an object of $\calC_U$
and $Q$ is an object of $\calD_U$, then there is a tensor product
object $C\ot D$ in $[\calC\ot\calD]_U$. We want any automorphism $g$ of $C$ or of $D$
to induce an automorphism of $C\ot D$, and these automorphisms should coincide.
However, such tensor
products of objects do not give enough objects to have the gluing properties
of a sheaf of
groupoids, so one needs to add more objects by gluing as follows. Let $(U_\a)$ be some open
covering of $U$ over which we have objects $P_\a
\in\calC_{U_\a}$ and $Q_\a\in\calD_{U_\a}$. Assume we have isomorphisms
$\phi_{\a\b}:[P_\b]_{/U_{\a\b}}\tilde{\to}[P_\a]_{/U_{\a\b}}$
and $\psi_{\a\b}:[Q_\b]_{/U_{\a\b}}\tilde{\to}[Q_\a]_{/U_{\a\b}}$.
We do not assume that the $\phi_{\a\b}$ and the $\psi_{\a\b}$ satisfy
the cocycle condition; instead we require that the defects of the cocycle
condition compensate each other, in other words if we have
$\phi_{\a\b}\ot\phi_{\b\g}\ot\phi_{\g\a}=c_{\a\b\g}$, we require that
$\psi_{\a\b}\ot\psi_{\b\g}\ot\psi_{\g\a}=c_{\a\b\g}^{-1}$. Then we postulate
that the objects $P_\a\ot Q_\a$ can be glued together, via the transition
isomorphisms $\phi_{\a\b}\ot\psi_{\a\b}$, to an object of $[\calC\ot\calD]_U$.
This gives enough objects in the category $[\calC\ot\calD]_U$ so that it is a DD gerbe.

The inverse $\calC^{\ot-1}$ of a gerbe $\calC$ has objects over $U$ given by the gerbe
equivalences from $\calC$ to the trivial gerbe. Thus there is an equivalence of gerbes
from $\calC^{\ot-1}\ot\calC$ to the trivial gerbe, given by evaluating an equivalence
of gerbes on an object of $\calC$.

The operations of tensor product of DD gerbes and inverse of a gerbe correspond
to the group structure on $H^2(M,\un\BC^*)$ (cf. Proposition 1).

It is harder to define the pull-back of a gerbe than the pull-back of gerbe data.
For  a smooth mapping $f:N\to M$ and for a gerbe $\calC$ on $M$, the pull-back gerbe
$f^*\calC$ is characterized by the fact that any object $P\in\calC_U$ yields an object
of $(f^*\calC)_{f^{-1}(U)}$. The construction of $f^*\calC$ involves several steps, and
we only give an outline. We can first use the procedure for pulling back a sheaf of
groups. Thus for
$V$ open in $N$, one can define the category $\calD_V$ to be the direct limit of the
categories
$\calC_U$, where $U$ runs over open neighborhoods of $f(V)\subseteq M$. There are
two problems with this construction.
First the automorphisms of objects are wrong: they are functions from open subsets
of $M$ to $\BC^*$, rather than from open subsets of $N$ to $\BC^*$. One remedies
this problem by formally enlarging the isomorphisms in $\calD$ so that the
automorphism group
of any  $P\in\calD_V$ is equal to the group of smooth functions $V\to\BC^*$.
Secondly there are not enough objects: one needs to add objects obtained by gluing
objects defined on an open covering $(V_\a)$ of $V\subset N$, together with gluing
data using the enlarged notion of isomorphisms. After this second step, one gets a
DD-gerbe on
$N$. We refer to [{\bf Br1, \S 5.2}] for details of the construction of the
pull-back.

The $0$ and $1$-connections take a different meaning in the context of gerbes
viewed as sheaves of groupoids. A $0$-connection becomes what we call a connective
structure in [{\bf Br 1}]: this associates to any object $P\in\calC_U$
a sheaf $Co(P)$ on which the complex-valued $1$-forms $A^1$ operate
``locally simply-transitively''. In more details, for any open set $U$, the set
$A^1(U)$ of
$1$-forms on $U$ operates on the sections of $Co(P)$ over $U$. For any $x\in M$,
there exists $U\ni x$ such that for any $x\in V\subseteq U$, $Co(P)(V)$ is isomorphic
to $A^1(V)$, where $A^1(V)$ acts on itself by translations. One says that
the sheaf $Co(P)$ is a torsor under the sheaf  $\un A^1_M$ of $1$-forms.
The trivial example is the trivial $0$-connection on the trivial gerbe. Then an object
over $U$ is a line bundle $L\to U$, and $Co(L)(U)$ is defined to be the set of connections
on the line bundle $L$.

For instance, if
$\calC$ is the trivial gerbe, then
$P$ is a line bundle
and we can take $Co(P)$ to be the set of connections on the line bundle. 
If we consider the gerbe attached to a principal $G$-bundle $Q\to M$ and a central
extension $\tilde G$ of $G$, then we get a connective structure once we fix a connection
on $Q$; then for an object $\tilde Q\to U$ of $\calC_U$, the sheaf $Co(P)$
is the sheaf of connections on $\tilde Q$ which lift the connection on $Q$

Then a $1$-connection becomes what is dubbed \und{curving} in [{\bf Br1}].
This associates to a section $\na$ of the sheaf $Co(P)$ some $2$-form $K(\na)$
(called the curvature of $\na)$, in such a way that
$$K(\na+\a)=K(\na)+d\a$$
for any $1$-form $\a$. Then the $3$-curvature becomes the $3$-form $\O$ such that
$\O=dK(\na)$.

The notion of curving lead to the interesting notion of flat object of a gerbe.
A flat object of $\calC$ over $U$ is an object $P$ of $\calC_U$ equipped with a 
section $\na$ of the sheaf $Co(P)$ such that $K(\na)=0$. Such a flat object can exist
only if the $3$-curvature $\O$ vanishes. If this is the case, flat objects exist
locally. From the point of view of gerbe data,  gerbe data
$(\La_{\a\b},\phi_{\a\b},\th_{\a\b\g})$ are flat if the line bundles $\La_{\a\b}$
are flat and the $\phi_{\a\b}$ and $\th_{\a\b\g}$ are horizontal.
Flat gerbes are classified by the cohomology group $H^2(M,\BC^*)$. Indeed if we pick
an open covering $(V_\a)$ and a flat object $(P_\a,\na_\a)$ over each $V_\a$,
as well as  isomorphisms $\psi_{\a\b}:P_\b\tilde{\to} P_\a$ which carry $\na_\b$
to $\na_\a$, then  $\th_{\a\b\g}=\psi_{\a\b}\psi_{\b\g}\psi_{\g\a}$ is an
automorphism of $P_\g$ which maps $\na_\g\in Co(P_\g)$ to itself; this implies
that $\th_{\a\b\g}$ is (locally) constant, hence yields a \v Cech cohomology
class in $H^2(M,\BC^*)$. If we fix the $3$-curvature $\O$ of a DD-gerbe, then the
DD-gerbes with $0$ and $1$-connection whose $3$-curvature is $\O$ differ from one another
by tensoring by a flat gerbe: hence they are parameterized by $H^2(M,\BC^*)$.

Curvings also occur in the geometric interpretation of the group-valued
moment maps of [{\bf A-M-M}] in terms of gerbes. This goes roughly as follows:
we suppose given on the compact simple Lie group $K$ a gerbe $\calC$ together with $0$
and $1$-connections, whose $3$-curvature is $2\pi i$ times the Chern-Simons
form. Given a smooth $K$-action on the manifold $M$ and an invariant $2$-form $\o$ on $M$, a
group-valued moment map is a $K$-equivariant mapping $\mu:M\to K$ (where
$K$ acts on itself by conjugation), together with an object $P$ of the pull-back gerbe
$\mu^*\calC$ and a section $\na\in Co(P)$ such that $K(\na)=\o$. There is an extra condition
in [{\bf A-M-M}, Def. 2.2], which has to do with the object $P$ and $\na$ being $K$-equivariant.

In this paper we will mostly deal with holomorphic gerbes over a complex manifold
$M$. This means, in terms of gerbe data, that the line
bundles $\La_{\a\b}$, the isomorphisms $\phi_{\a\b}$ and the sections $\th_{\a\b\g}$
are holomorphic. In terms of gerbes, a holomorphic gerbe is a sheaf of groupoids
$\calC$ such that the group of automorphisms of an object $P$ of $\calC_U$ is the
group of holomorphic functions $U\to \BC^*$. A holomorphic gerbe data yields a
standard smooth data by simply forgetting that the data 
$(\la_{\a\b},\phi_{\a\b},\th_{\a\b\g})$ are. From the point of view
of sheaves of groupoids, a holomorphic gerbe yields a smooth gerbe
by first extending the isomorphisms between objects to accept smooth functions
to $\BC^*$ as automorphisms, and
then adding new objects  obtained by gluing with the help of the new automorphisms.
Similarly we have the notion of a holomorphic $0$-connection, meaning that the
connections $D_{\a\b}$ on the $\La_{\a\b}$ are holomorphic. On the side of sheaves of
groupoids, one needs to attach to each object $P$ of any $\calC_U$ a torsor
$Co_{hol}(P)$ under the sheaf $\O^1_M$ of holomorphic $1$-forms.
The one recovers the torsor  under the sheaf $\un A^1_M$ of 
smooth $1$-forms by enlarging
the torsor $Co_{hol}(P)$ accordingly. For $1$-connections to be holomorphic, one
requires the corresponding $2$-forms to be holomorphic.

The classification now involves \v Cech cocycles $g_{\a\b\g}$ which are holomorphic.
In other words, they are \v Cech $2$-cocycles with values in the sheaf $\calO_M^*$
of invertible holomorphic functions. Thus equivalence classes of holomorphic gerbes are
classified by the
\v Cech cohomology group $H^2(M,\calO^*)$. Similarly, holomorphic gerbes with
holomorphic $0$-connection are classified by the \v Cech hypercohomology group
$H^2(M,\calO^*\mapright{d\log}\O^1_M)$, and holomorphic gerbes with holomorphic $0$
and $1$-connections by the 
\v Cech hypercohomology group
$H^2(M,\calO^*\mapright{d\log}\O^1_M\mapright{d}\O^2_M)$.
These are (holomorphic) Deligne cohomology groups; see [{\bf Br1}] [{\bf Br2}]
for details.

\sec{2. The Grothendieck manifold}

Let $G$ be a connected reductive algebraic complex Lie group.
Let
$X$ be the the variety of Borel subgroups of $G$. Recall that
the $G$-orbits in $X\times X$ are canonically parameterized
by a finite group $W$ called the Weyl group of $G$.
Denote by $Y_w$ the orbit corresponding to $w\in W$. This a
version
of the Bruhat decomposition. Thus given two Borel subgroups
$B,B^{\pr}$ of $G$, we say that $(B,B^{\pr})$ are in position
$w$ (notation $B\mapright{w}B'$) if $(B,B^{\pr})\in
Y_w$. It follows that $W$ acts naturally on the set $X^T$ of
Borel subgroups containing a given Cartan subgroup $T$.
Given $B\supset T$ and given $w\in W$, there exists a unique
$B^{\pr}
\supset T$ such that $B\mapright{w}B'$. We let
$w\cdot B=B^{\pr}$. This action is simply-transitive.
Now it easily follows that given $T\subset B$ we can identify
the geometric Weyl group $W$ with the concrete Weyl group
$W_T=N(T)/T$. Indeed both groups act on $X^T$ and the two
actions
commute and are simply-transitive. Thus the choice of a point
$B\in X^T$ yields an isomorphism between the two groups.

Now let
$\tilde G\subset G\times X$ be the closed subset consisting of
pairs
$(B,g)$ where $g\in B$. Then the projection
$\tilde G\to X$ is a locally trivial fiber bundle with fiber
over $B\in X$ equal to $B$. Hence $\tilde G$ is an algebraic
bundle of groups over $X$ and is a smooth algebraic variety.
We call $\tilde G$ the \und{Grothendieck manifold}; it was introduced by
Grothendieck for the purpose of resolving simultaneously
the singularities of all closures of conjugation orbits in $G$.
On the other hand, the projection $q:\tilde G\to G$ is a
proper algebraic mapping. Let $G^{reg}\subset G$ be the open
set of regular semisimple elements.

\lem{1}{The restriction of $q$ to $G^{reg}$ is a finite Galois
cover with Galois group $W$.}

\proof{There is a natural action of $W$ on $q^{-1}(G^{reg})$:
given $(g,B)\in q^{-1}(G^{reg})$, there is  a unique Cartan
subgroup $T$ containing $g$. The fiber of $\tilde G$ over $g$
thus identifies with $X^T$, which admits a natural action of
$W$ as we discussed above. It is clear that this action
is simply-transitive on the fibers of $q$, and makes
the restriction of $q$ to $G^{reg}$ into a Galois covering. \qed}

The action of $W$ on the open set $q^{-1}(G^{reg})$ of
$\tilde G$ does not extend to a global action on $\tilde G$.
However there is a natural action of $W$ on the cohomology
of $\tilde G$. This can be seen as follows.
We have a projection
map $p:\tilde G\to X\times T$, which is a
homotopy equivalence. Now there is a well-known action of $W$
on $X=G/B$, which depends on the choice of  maximal compact
subgroup $K$ of $G$. Indeed, introducing the maximal torus
$T_c=T\cap K$ of $K$, we have $X\simeq K/T_c$ and $W_T$
is the Weyl group of $T_c$, which acts on $K/T_c$ by right
multiplication.

Then we have the diagonal action of $W$ on the product $X\times T$.
Thus we get an action of $W$ on $H^*(\tilde G)=H^*(X\times T)$.
The restriction of $p$ to $q^{-1}(G^{reg})$ is $W$-equivariant.

The submanifold $(K/T_c)\times T_c$ of $\tilde G$ has half the real dimension
of $\tilde G$ and is a $W$-equivariant deformation retract of $\tilde G$. We will call it
the
\und{core} of $\tilde G$.

In the following Proposition, cohomology groups are taken
with coefficients in a field of characteristic $0$.

\prop{3}{The pull-back map induces an isomorphism
of $H^*(G)$ with $H^*(\tilde G)^W$.}

\proof{There are several ways of proving this result.
The most concrete is to first observe that since $q$
is a proper map which is generically finite, the
pull-back map $H^*(G)\to H^*(\tilde G)$ is injective.
Clearly its image is contained in the $W$-invariant subspace.
Now it is easy to compute $H^*(\tilde G)$ together with its
$W$-action. For this we use the core $X\times T_c$ of $\tilde G$.
 We have 
$H^*(X\times T_c)=H^*(K/T_c)\ot H^*(T_c)$. We have natural actions of $W$
on $H^*(K/T_c)$ and on $H^*(T_c)$. The $W$-module $H^*(K/T_c)$ is
isomorphic to the regular representation. Hence the
$W$-invariant subspace of
$H^*(X)\ot H^*(T_c)$ identifies with $H^*(T_c)$.
Since $H^*(G)$ and $H^*(T_c)$ have the same dimension,
the statement follows. \qed}

\cor{2}{The pull-back map induces an isomorphism
of the equivariant cohomology ring $H^*_G(G)$ with
$H^*_G(\tilde G)^W$.}

\proof{This follows immediately from Proposition 1,
using the spectral sequence from ordinary to equivariant
cohomology. \qed}

\rm We can compute $H^*_G(\tilde G)^W$ as follows.
We have 
$$H^*_G(\tilde G)=H^*_K(\tilde G)=H^*_K(X)\ot H^*(T_c)
=R(T_c)\ot H^*(T_c)=\BQ[X^*(T_c)]\ot \wedge^*X^*(T_c),$$
where $X^*(T_c)=Hom(T_c,S^1)$ is the character group of $T_c$
and $R(T_c)$ is the representation ring of $T_c$.
We can identify $X^*(T_c)$ with the algebraic character group $X^*(T)$
of the algebraic torus $T$ and $R(T_c)$ with the algebraic representation ring $R(T)$.
The tensor product algebra $\BQ[X^*(T_c)]\ot \wedge^*X^*(T_c)$ identifies with the ring
$\Omega^*(R(T)\ot\BQ)$ of Grothendieck differentials of the algebra
$\BQ[X^*(T)]=R(T)\ot\BQ$. The following result is proved in
[{\bf B-Z}].

\lem{2}{We have:
$$\Omega^*(R(T)\ot\BQ)^W=\Omega^*(R(T)^W\ot\BQ).$$ }
Thus we obtain

\prop{4}{We have a natural isomorphism of graded algebras:
$$H^*_G(G)=\Omega^*(R(T)^W\ot\BQ).$$}

The nice thing about the core
$(K/T_c)\times T_c$  of $\tilde G$ is that its equivariant 
cohomology can easily be represented by explicit equivariant
differential forms. The cohomology of $X$ is generated by
Chern classes
of line bundles.  Then any character $\chi\in T$ extends
uniquely to a character of $B$ hence defines a 
$G$-homogeneous holomorphic line bundle $L(\chi)$
over $G/B$. We can get a hermitian structure
on $L(\chi)$ by picking 
a maximal compact
subgroup $K$ of $G$. Then for any $\chi\in
X(T)$ the line bundle $L(\chi)$
acquires a hermitian structure  from
the facts that $X\simeq K/T_c$  and the restriction of
$\chi$ to $S\subset T_{sc}$ takes values in the circle group
$\BT\subset\BC^*$. 

Since the line bundle $L(\chi)$ over $X$ is homomorphic and
hermitian, it has a canonical connection. Let
$R_\chi$ be the curvature of this connection. This is a
purely imaginary $K$-invariant closed
$2$-form, and its value at the base point can be computed as
follows. The tangent space to $K/T_c$ is the quotient
space
$\fk/\ft_c$; for $\xi\in\fk$ we denote by $\bar\xi$
the corresponding tangent vector. Then  for elements
$\xi,\eta$ of
$\fk$ we have:

$$R_\chi(\bar\xi,\bar\eta)=d\chi([\xi,\eta]).$$

The $2$-form $\o_{\chi}={1\over 2\pi\sqrt{-1}}R_{\chi}$
has integer periods, so its cohomology class belongs
to $H^2(X,\BZ)\subset H^2(X,\BC)$. We note that the curvature
$R$ of the principal bundle $K\to K/T_c$ is the $\ft$-valued $2$-form
$$R(\bar\xi,\bar\eta)=p_\ft\chi([\xi,\eta])$$
where $p_\ft:\fk\to\ft$ is the projection; then we have
$R_{\chi}=\lan d\chi,R\ran$.

The $2$-form $\o_{\chi}$ can be extended to a closed
$K$-equivariant $2$-form on $X$, namely
$$(\o_{\chi},\mu_{\chi})$$
where $\mu_{\chi}:X=K/T_c\to\fk^*$ is the moment map
$$\mu_{\chi}(kT_c/T_c)={1\over 2\pi\sqrt{-1}}\,Ad^*(k)d\chi.$$
Here $d\chi\in \ft_c^*$ is extended to an element of
$\fk^*$ which vanishes on the orthogonal of $\ft_c$.

So $(\o_{\chi},\mu_{\chi})$
yields an equivariant differential form, whose class
in $H^2_K(X)=H^2_G(X)$ is the opposite of 
the equivariant first Chern class $c_1^K(L(\chi))$.

Now any character $u$ of $T$ gives rise to the equivariant
$1$-form ${1\over 2\pi\sqrt{-1}}{du\over u}$ over $T$, which represents the
 cohomology class of $u$ in $H^1(T,\BZ)=X(T)$.

Then, denoting by $A^*_K(X\times T)=[A^*(X\times T)\ot
S(\fk)^*]^K$ the complex of
$K$-equivariant
differential forms, we have an algebra map
$${\bf c}:\BZ[X(T)]\ot\wedge^* X(T)\to A^*_K(X\times T)$$
such that
$${\bf c}(\chi\ot 1)=p_1^*(\o_{\chi},\mu_{\chi})\,\,,\,\,
{\bf c}(1\ot \chi)={1\over 2\pi \sqrt{-1}}p_2^*{d\chi\over\chi},$$
where $p_1:X\times T\to T$ and $p_2:X\times T\to T$ are the two projections.

It is interesting to compare the canonical equivariant
$3$-forms on $G$ and on $\tilde G$ attached to
a $W$-invariant bilinear form $b$
on the cocharacter group $X_*(T)=X^*(T)^*$. Note that $X_*(T)$ is the group of algebraic
homomorphisms
$\la:\BC^*\to T$. Then
$b$ extends
uniquely to an invariant complex-valued bilinear form $b$ on $\fg$.

On one hand, there is the well-known $G$-equivariant
$3$-form
$$(\nu,\a).$$
To describe it, introduce the left (resp right) invariant
Maurer-Cartan forms  $\th$ (resp. $\bar\th$) over $G$.
 Then
$\nu$ is the Chern-Simons
$3$-form
$$\nu={1\over 12}b(\th,[\th,\th])$$
and $\a$ is the $\fg$-valued $1$-form
such that
$$\a={1\over 2}(\th+\bar\th).$$
We will consider the restriction of $(\nu,\a)$ to a
$K$-equivariant differential form on $G$.
This $K$-equivariant differential form is closed.

On the other hand, we can view $b$ as an element
of $X^*(T)\ot X^*(T)\subset \BZ[X^*(T)]\ot\wedge^*
X^*(T)$ and then we have the $K$-equivariant differential form
${\bf c}(b)$ on $X\times T$.

By pulling back both equivariant $3$-forms to $\tilde X$ we
can compare them:

\prop{5}{We have the equality of $K$-equivariant differential
forms on $\tilde X$:
$$q^*(\nu,\a)-{1\over 2}p^*{\bf c}(b)=d_K(p^*\b)$$
where $\b$ is the invariant $2$-form on $ X\times T$ defined as follows.
$\b(v,w)$ vanishes if $v$ or $w$ is tangent to the factor $T$.
For $\xi, \eta\in\fk$ we have 
$$\b_{1,t}(\bar\xi,\bar\eta)={1\over 2}\left(
b(\xi,Ad^t(\eta))-b(\eta,Ad^t(\xi)\right).$$}

This follows from Theorem 6.2 in [{\bf G-H-J-W}] or Proposition 3.1 in [{\bf A-M-M}].

We note a variant of these constructions.
To get the full set of homogeneous line
bundles, it is important to consider line bundles
over $X$ which are equivariant under more general groups.
Thus let $f:H\to G$ be an algebraic homomorphism of complex algebraic Lie
groups. We will assume that $Ker(f)$ is central in $H$ and that
$G=Z(G)\cdot Im(f)$, where $Z(G)$ is the center of $G$. Then we have $X\simeq H/f^{-1}(B)$.
Thus any character
$\chi$ of $f^{-1}(B)$ defines an equivariant line bundle
$L(\chi)$ over $X$. Now $f^{-1}(B)$ has the same group of characters
as its subgroup $S=f^{-1}(T)$. If we pick a maximal compact subgroup $L'$
of $f^{-1}(K)$, the product $L=L'\cdot Z(H)$  will also be
maximal compact in
$H$. As the restriction of $\chi$ to $S\cap L$ will take  values
in $\BT$, it follows that $L(\chi)$ has a hermitian structure.
Hence we can extend the construction of the $2$-forms $R_{\chi}$
and $\o_{\chi}$ to this case.
We thus obtain an algebra map 
$$S^*(X^*(S))\ot \wedge^* X(T)\to A^*_L(X\times T).$$

A very interesting case is the following.
Introduce the derived subgroup
$G^{\pr}$ of $G$ and its simply-connected covering $G_{sc}$,
which is a complex semisimple algebraic group.
We then have $X\simeq G_{sc}/B_{sc}$, where $B_{sc}=f^{-1}(B)$.
Then $S$ be a maximal torus of $G_{sc}$ contained in
$B_{sc}$. Then any character $\chi\in S$ extends
uniquely to a character of $B_{sc}$ hence defines a 
$G_{sc}$-homogeneous holomorphic line bundle $L(\chi)$
over $G_{sc}/B_{sc}$. The $G_{sc}$-equivariant cohomology of
$\tilde X$ and of $X\times T$ is then equal
to $$S^*(X^*(T_{sc}))\ot \wedge^* X(T).$$
This algebra maps to the algebra of $K_{sc}$-equivariant
differential forms on $X\times T$.

\sec{3. Construction of gerbes over $\tilde G$}

As in section 1, let $f:H\to G$ be an algebraic homomorphism of complex
algebraic Lie groups which is an isomorphism modulo centers. We will give a general
construction for
$H$-equivariant holomorphic gerbes on $\tilde X$. The basic data we will use is
an element $b\in X^*(S)\ot X^*(T)$, which we often view as a bilinear form
$$b:X_*(S)\ot X_*(T)\to \BZ.$$
We will associate to $b$ the corresponding linear map
$\un b:X_*(T)\to X^*(S)$. Thus for any $\la\in X_*(T)$
we have the $H$-equivariant line bundle $L(\un b(\chi))$ on $X$.

We will first give two constructions of  gerbe data on $\tilde X$.
Given
a line bundle $L$ over a manifold $X$ and a smooth
function $f:X\to\BC^*$, there is a gerbe attached to $L$ and $f$, which was discussed
and exploited in [{\bf B-M2}]. This gerbe is given by a cup-product construction. We will
discuss two methods to construct gerbe data corresponding to $L$ and $f$. In the first
method, we pick an open covering
$(V_\a)$ of
$X$ over which we have a branch
$\log_\a~f$ of a logarithm of
$f$. Then we have $\log_\b ~f-\log_\a~f=2\pi i m_{\a\b}$ for $m_{\a\b}\in\BZ$.
Then we set $\La_{\a\b}=L^{\ot m_{\a\b}}$ with the obvious choice of $\phi_{\a\b}$
(see \S 1). The
trivialization
$\th_{\a\b\g}$ is obvious from the cocycle property
$m_{\b\g}+m_{\a\g}+m_{\a\b}=0$.

The second method uses the Deligne line bundle $(g,f)$ attached to
two invertible functions $f,g$ [{\bf De}]. To describe $(g,f)$, it is enough to work in
the universal case where $X=\BC^*\times\BC^*$, and $g=x, f=y$.
Then we start with the trivial bundle over   the
universal covering space $\BC^*\times\BC$ of $\BC^*\times\BC^*$ with covering map
  $\pi:=Id\times exp:\BC^*\times\BC \to
\BC^*\times\BC^*$. The group of deck transformations is  $\BZ$, where
the generator $T$ acts on $(x,w)$ by $(x,w)\mapsto (x,w+2\pi i)$. We let
$T$ act on the trivial line bundle over $\BC^*\times\BC$ by multiplication
by the invertible function $x^{-1}$. In other words, a section of $(f,g)$
over some open set $U$ is a function $h:\pi^{-1}(U)\to \BC$ such that
$h(x,w+2\pi i)=x^{-1}h(x,w)$. The Deligne line bundle is equipped with a connection
$\na$ which is characterized by the equation
$\na(h)=dh+(2\pi i)^{-1}wx^{-1}h$. The curvature is
$$(2\pi i)^{-1}{dx\over x}\wedge {dy\over y}.$$
Note that by construction, a (local) logarithm $\log f$ determines
a non-vanishing section of $(g,f)$, which we will denote by $(g,\log f]$. We have the
rule $(g,\log f+2\pi im]=x^m(g,\log f]$. We refer to [{\bf De}] for more information
on the Deligne line bundle, which Deligne invented in relation with algebraic K-theory
and regulator maps. For our purposes (as in [{\bf Br1}] [{\bf B-M1}]), we view it as a
very nice way
to quantize the $2$-torus or its complexification.

Returning to our data $(X,L,f)$, we give our second construction
of gerbe data. We pick an open covering
$(W_{\a})$ over which $L$ has a non-vanishing section $s_\a$.
Then we have $s_\b=s_\a g_{\a\b}$ for the transition cocycle
$g_{\a\b}$. We then define the line bundle $\La_{\a\b}=(g_{\a\b},f)$
with obvious choices of $\th_{\a\b\g}$. Each $\La_{\a\b}$ is equipped with a
connection, for which $\th_{\a\b\g}$ is horizontal. Thus we have a
$0$-connection on the gerbe data.
The equivalence of this gerbe with the previous one is easy to prove:
it amounts to showing that the gerbe data given by the line bundles
$(g_{\a\b},f)\ot L^{\ot m_{\b\a}}$ is trivial. This can be done by
using the trivializing sections $(g_{\a\b},\log_\a f]\ot s_\a^{\ot m_{\b\a}}$.

We will denote by $(L,f)$ the gerbe we have constructed.

Now in our situation we have a tensor $b$ in $X^*(S)\ot X(T)$.
If we write $b=\sum_b\chi_j\ot \ze_j$, then for each $j$ we have
the $H$-equivariant line bundle $L(\chi_j)$ over $X$ and the invertible
function $\ze_j$ on $T$, hence we have the tensor product gerbe $\ot_j(L(\chi_j),\ze_j)$
on the product $X\times T$. Furthermore this gerbe is $H$-equivariant
(cf. the appendix for this notion). We can take its pull-back under the
map $p:\tilde G\to X\times T$ to get an $H$-equivariant gerbe
with $0$-connection on $\tilde G$. We will denote this gerbe by
$\tcalC=\tcalC_b$.

We need however to make sure that this construction does not depend on the
expression of $b$ as a sum of tensors. This is true up to natural
isomorphisms. For instance, using the second construction, the main point is
that for a complex manifold $Y$ and for an
element $\g\in Hol(Y)^*\ot Hol(Y)^*$, there is a canonical Deligne line
bundle $(\g)$ which is isomorphic to $\ot_j (g_j,f_j)$ for any
expression of $\g$ as $\sum_j g_j\ot f_j$.

\prop{6}{If the bilinear form $b$ is $W$-invariant, then
the restriction of
the gerbe
$\calC_b$ to
$\tilde G^{reg}$ is $W$-equivariant.}

\proof{For this purpose, it is best to use the description
of $\tilde G^{reg}$ as $G\times^T T^{reg}=(G/T)\times T^{reg}$. Then the
$W$-action on this contracted product can  be described as a diagonal action.
It is clear that with respect to the action of $W$ on $G/T$, we have a
canonical isomorphism
$w^*L(\chi)\tilde{\to}L(w^{-1}\chi)$. On the other hand, the $W$-action
on $T$ transforms the characters of $T$ according to the action of $W$
on $X(T)$. The result then follows since the construction of $\tcalC$
does not depend on the expression of $b$ as a sum of tensors. \qed}

Our gerbe $\calC_b$ is such that the sum $b_1+b_2$ ot two tensors
leads to the tensor product gerbe $\calC_{b_1}\ot\calC_{b_2}$.
The gerbe $\calC_b$ has a simple behavior under the inverse map
$\iota:\tilde G\to\tilde G$ given by $\iota(g,B)=(g^{-1},B)$:

\lem{3}{The pull-back $\iota^*\calC_b$ is equivalent to $\calC_b^{\ot-1}=\calC_{-b}$.}

\proof{We have the commutative square
$$\matrix{G&\mapright{\iota}&G\cr
\mapdown{}&&\mapdown{}\cr
T&\mapright{inv}&T}$$
where $inv:T\to T$ is the inverse map.
On the other hand, $\iota$ is compatible with the projection $\tilde G\to X$.
Then using the construction of $\calC_b$ as $\ot_j {L(\chi_j,\ze_j)}$
we see that $\iota^*\calC_b$ is $\ot_j {L(\chi_j,inv^*\ze_j)}=\ot_j
{L(\chi_j,-\ze_j)}=\calC_b^{\ot-1}$. \qed}

We conclude with two reasonably concrete description of our gerbe $\tcalC$.
The exponential map $exp:\ft\to T$ is a covering whose 
Galois group is $X_*(T)$; for $\la\in X_*(T)$
we denote by $t_{\la}$ the corresponding  deck transformation $\ft\to\ft$ which is
translation by
$2\pi i[d\la(1)]$. Denote by
$\kappa:\tilde G\to T$  and by $p_2:\tilde G\to X=G/B$
the projection maps. Let $U\subset \tilde G$ be an open set
and let $\hat U=U\times_T\ft\to U$ be the pull-back covering of $U$. For $\la\in X_*(T)$ we
still have the translation automorphism $t_{\la}$ of $\hat U$. Then an object of
$\tcalC(U)$ is a holomorphic line bundle $\calL$ over $\tilde U$ together
with isomorphisms
$$\b_{\la}:t_{\la}^*\calL\tilde{\to}\calL\ot L(-b(\la))$$
for $\la\in X_*(T)$, which satisfy the transitivity condition
$$\b_{\la_1+\la_2}=(t_{\la_1}^*\b_{\la_2})\b_{\la_1}.$$
Here we still denoted by $L(-b(\la))$ the pull-back to $\hat U$
of the corresponding line bundle over $X$.

The second description is derived from the first, but it
directly gives gerbe data for
an open covering of $\tilde G$. We cover
$T$ by open sets
$V_\a$ over which we have a section $\log_\a:V_\a\to \ft$ of the exponential map
$exp:\ft\to T$. Over $V_{\a\b}$ we have $log_\b-\log_\a=2\pi i \la_{\a\b}$
for some $\la_{\a\b}\in X_*(T)$. Then we cover
$\tilde G$ by the open sets $U_\a=\kappa^{-1}(V_\a)$. Over $U_{\a\b}$ we consider
the line bundle $\La_{\a\b}=p_2^*L(b(\la_{\a\b})$. With obvious
choices for $\th_{\a\b\g}$ this gives gerbe data.

This last description is also useful to give a $1$-connection on the gerbe
restricted to the core $(K/T_c)\times T_c$ of $\tilde G$.
Over the intersection of the core with $U_\a$ we have the $2$-form
$F_\a=b(p_2^*\log_\a,p_1^*R)$. The $3$-curvature is the $3$-form described in \S 2
as the equivariantly closed $3$-form ${\bf c}(B)$.
We conjecture that if we add $2\b$ to this $1$-connection, the resulting $1$-connection
extends holomorphically to $\tilde G$. This should be expected, as the corresponding
$3$-form is equal
to $2\pi i\nu$ by Proposition 5, and this extends to a holomorphic $3$-form over $G$.

\sec{4. Descent to $G$: the case of $SL(2)$.}

We have constructed a gerbe $\tcalC$ with $0$-connection on $\tilde G$,
and we now want to descend it to $G$. The first step is to construct
the gerbe $\calC$ over the open set $G^{reg}$. There we can use Lemma 1 which
gives us the Galois covering $\tilde G^{reg}\to G^{reg}$ with Galois
group $W$. Since the gerbe $\tcalC$ on $\tilde G^{reg}$ is $W$-equivariant,
it automatically descends to a gerbe on $G^{reg}$. This is much easier to
describe for gerbes than for gerbe data. An object of the gerbe $\calC$ over
an open set $U\subset G^{reg}$ will be an object $P$ of $\tcalC$ over
$q^{-1}(U)$ which is $W$-equivariant, that is equipped with isomorphisms
$\eta_w:w^*P\tilde{\to}P$ such that $\eta_{w_1w_2}=w_1(\eta_{w_2})\eta_{w_1}$
(see the Appendix for a discussion of equivariant gerbes and equivariant
objects).
Then we have:

\lem{4}{This construction describes a gerbe $\calC$ on $G^{reg}$.}

Next we want to extend this gerbe to $G$. First we examine the case
where $G$ is a complex connected semisimple algebraic group of dimension
$3$, that is $G=SL(2,\BC)$ or $G=PGL(2,\BC)$. Let $Z(G)$ denote the center of
$G$, which is a finite group of order $2$ or $1$. We have  the following
nested open subsets of $G$:

$$G^{reg}\subset V\subset G,$$
where $V=G\setminus Z(G)$. Over $V$, the mapping $\tilde V:=q^{-1}(V)\to V$
is a ramified double covering with Galois group $W=\BZ/2=\{ 1,\tau\}$.
So we are led to the question of descending the $W$-equivariant gerbe $\tcalC$
on $\tilde V$ to a gerbe on $V$. The complement $Y$ of $G^{reg}$ in $V$
is a smooth hypersurface: for $SL(2)$ it has $2$ components $Y_1,Y_2$
corresponding to matrices all of whose eigenvalues are $1$ resp. $-1$,
and for $PGL(2)$ it is connected.

Now consider a possibly ramified double covering $p:\tilde V\to V$ with involution $\tau$,
and a $\BZ/2$-equivariant holomorphic gerbe $\tcalC$ over $\tilde V$. First
of all, if the covering is not ramified, we construct a holomorphic gerbe
$\calC$ over $V$ whose objects over $U\subset V$ are the $\BZ/2$-equivariant
objects of $\tcalC_{p^{-1}(U)}$ (see Appendix for equivariant objects); the morphisms
are $\BZ/2$-equivariant isomorphisms. If $\tcalC$ has an equivariant connective structure
$P\mapsto Co(P)$,
then for any equivariant object $P$ of $\tcalC_{p^{-1}(U)}$ the $\O^1_{p^{-1}(U)}$-torsor
$Co(P)$ is $\BZ/2$-equivariant, hence it descends to an $\O^1$-torsor over $U$. Thus
$\calC$ acquires a connective structure.

Now consider the case where the covering $p:\tilde V\to V$ is ramified
along a smooth hypersurface $Y$. In that case there is an obstruction to descending 
the  equivariant gerbe $\tcalC$: this consists of an element of $\BZ/2$ attached
to each component $Y_j$ of $Y$. The description of this integer mod $2$ is
purely local along $Y_j$; thus we may assume the gerbe is trivial, so is the
gerbe whose objects are line bundles. Then the action of $\tau$ on $\tcalC$
must given by $\tau(L)=L\ot \La$ for some line bundle $\La$
on the complement of $Y_j$ (cf. Proposition 2). The
constraint that $\tau^2$ should be (isomorphic to) the identity means
that there is an isomorphism $\phi:\La\ot\tau^*\La\tilde{\to} \bf 1$.
Now the obstruction arises when we look for a local section $s$ of $\La$
around some point of $Y_j$
such that $\phi(s\ot\tau^*(s))=1$. Indeed, first take any holomorphic
section $\s$ of $\La$ and consider the order $d$ of
$f:=\phi(\s\ot\tau^*(\s))$ along $Y_j$. When we multiply $\s$ by
a meromorphic function $g$ (with possible pole along $Y_j$), we change this order into
$d+2l$, where
$l$ is the order of
$g$ along $Y_j$. Thus the residue modulo $2$ of $d$ is an intrinsic
invariant, and is our obstruction. It can clearly be measured by
restricting the whole geometric situation (including the gerbe)
to some small $\tau$-invariant disc which meets $Y_j$ transversally at the
origin.

If the obstruction vanishes, then we have a holomorphic gerbe over $V$, whose objects
over $U$ are again the $\BZ/2$-equivariant objects of $\tcalC$ over $p^{-1}(U)$.
Then an equivariant connective structure on $\tcalC$ will induce one on $\calC$ just as
in the unramified case.

When $G=SL(2,\BC)$, there are $2$ divisors $Y_1,Y_2$ to consider. The open set 
$G^{reg}$ is isomorphic to $(G/T)\times T^{reg}$, where $T=\BC^*$
and $T^{reg}=\BC^*\setminus \{\pm 1\}$. The variable in $T$ will be denoted
by $z$. $W$ acts on $T$ by $z\mapsto z^{-1}$. The character group $X^*(T)$ is
generated by the identity character $\chi$. 
$b$ is determined by the integer $m=b(\chi\ot\chi)$. The line bundle $\La$ occurs
because to write down the gerbe data, we need to fix a local branch of the
logarithm of
$z$, and this will not in general be $W$-equivariant. Near $z=1$ we can make a
$W$-invariant choice of $\log(z)$, such that $|\log(z)|<\pi$, so the line
bundle $\La$ is trivial. However near $z=-1$ (corresponding to the divisor
$Y_2$), if we choose  the branch so that $\log(-1)=\pi i$, then we have
$\log(z^{-1})=-\log(z)+2\pi i$. It follows that the line bundle
$\La$ over a neighborhood of $Y_2$ in $G^{reg}$ is the pull-back
of the line bundle $L(m\chi)$ over $G/T$. We must then analyze the 
element of $\BZ/2$ attached to this line bundle.

The surface $G/T$ is an affine quadric, and as such it has two rulings
which are exchanged by $W$. We call the lines of these rulings lines of the first resp.
second kind.  The line bundle $\La$ is constant along the  lines of the first kind,
which are the fibers of the projection $G/T\to G/B$. Along the lines of the
second kind, the line bundle corresponds to the divisor $m[p]$ for a point $p$. Its
transform
$\tau^*\La$ is constant along the lines of the
second kind, and its restriction to lines of the first kind is attached to some point.
Now as we approach the divisor $Y_2$ along a transverse disc, we stay in a small
neighborhood of some line of the second kind (this is because the projection
to $G/B$ of our point in $G/T$ has a limit). It follows that the pull-back of $\La$
to this disc has a zero of order $m$ at the origin, while $\tau^*\La$ pulls back
to an equivariantly trivial bundle. Thus we conclude

\lem{5}{For $G=SL(2,\BC)$, the gerbe $\tcalC$ can descended from $\tilde V$
to $V=G\setminus \{ \pm 1\}$ iff $m=b(\chi\ot\chi)$ is even.}

The case of $PGL(2,\BC)$ is  different as there is only one divisor and the
corresponding obstruction vanishes automatically. This can be viewed in the following
way: the integer $m$ is even because $\display\chi\over 2$ is a character of $T$.

Next we need to extend our gerbe from $V$ to $G$.  We use
the following Hartogs type theorem from [{\bf Br2}].

\lem{6}{Let $X$ be a complex manifold, and let $Z$ be a closed complex subvariety
of codimension $\geq 3$. Then the restriction map
from holomorphic gerbes with $0$-connection on $X$ to those on $X\setminus Z$ is
an equivalence of categories.}

So we obtain a gerbe $\calC$ on $G$ equipped with a $0$-connection.

\sec{5. Descent to $G$: the general case.}  

We will first construct an extension of the gerbe $\calC$ over $G^{reg}$
to $G\setminus Z$, where $Z$ is a Zariski closed subset of codimension $\geq 2$.

The complement $Y$ of $G^{reg}$ in $G$ is a divisor, whose components
$Y_\a$ are indexed by the roots $\a$ up to the $W$-action; so if the semisimple part
of $G$ is simple, there is one component in the simply-laced case
and two otherwise. A general point of $Y_\a$ is $G$-conjugate to an
element $g$ with Jordan decomposition $g=su$, where

- $s$ is an element of $T$ such that $exp(\pm\a)(s)=1$ but no other root
is trivial on $s$;

- $u$ is a general unipotent element of the centralizer $Z_G(s)$.

Note that $Z_G(s)$ is (up to a finite group) the product of a torus
of dimension $r-1$ with a $3$-dimensional simple group $R$.

The $W$-covering $\tilde G^{reg}\to G^{reg}$ has an ordinary ramification
of order $2$ along each component $Y_\a$; the corresponding ramification
subgroup is the $\BZ/2$-subgroup generated by $s_\a$.

Now we have the $W$-equivariant gerbe $\tcalC$ and we wish to descend it
to an open subset of $G$ which at least meets each component $Y_\a$.
We studied in \S  3 the obstruction to doing this: it is an element
of $\BZ/2$ attached to each component $Y_\a$.  Now pick a general point
$g=su$ of $Y_\a$ as above, so we have a closed embedding $R\mapright{j} G$
where $j(h)=sh$, and the trace of the divisor $Y_\a$ on $R$ is a component
$V$ of the complement of $R^{reg}$ in $R$. We can lift $j$ to an algebraic
mapping $\tilde j:\tilde R\to\tilde G$ of Grothendieck manifolds as follows.
Fix a Borel subgroup $ B$ of $G$ containing $s$ and let $P_\a$ be the minimal
parabolic subgroup containing $b$ corresponding to the root $\a$.
Given a Borel subgroup $C$ in $R$, there is a unique Borel subgroup $B'$
of $G$ such that $C\subset B'\subset P_\a$. Then we set
$\tilde j(h,C)=(sh,B')$. Now we wish to describe the pull-back of the gerbe $\tcalC$
under this mapping $\tilde j$. For this purpose, we define a  maximal torus $T_R$
of $T$ as the intersection of $T$ with $R$. Then we define the algebraic group
$L=R\times_GH$ which maps to $R$, and we define a maximal torus $S_R$ of $L$ to
fit in the cartesian diagram

$$\matrix{L=R\times_GH&\to&H&\to&G\cr   
\mapup{}&&\mapup{}&&\mapup{}\cr
S_R&\to&S&\to&T}$$

Then the algebraic group homomorphism $L\to H$ satisfies our assumptions
and we can view the restriction $j^*\tcalC$ of $\tcalC$ to $R$ as
an $H_R$-equivariant gerbe.

We can then consider the restriction map
$X^*(S)\ot X^*(T)\to X^*(S_R)\ot X^*(T_R)$.
Denote by $b_R$ the image of $b\in X^*(S)\ot X^*(T)$ under this map. 

\lem{6}{The $L$-equivariant gerbe $j^*\tcalC$ over $R$ is the gerbe associated to the
element $b_R$ of $X^*(S_R)\ot X^*(T_R)$.}

It follows then that for each component $Y_\a$ of $Y$, the obstruction
in $\BZ/2$ can also be calculated in terms as the obstruction to extending the gerbe
on $R^{reg}$ along the divisor $j^*Y_\a$. We know from \S 3 that this obstruction
vanishes if $b_R(\check\a,\check\a)$ is even. Now this is this equal to
$b(\check\a,\check\a)$. Hence if we make the assumption

$$b(\check\a,\check\a)~{\rm is~even~ for~ any~ root}~\a\eqno(EV)$$
then we can descend $\tcalC$ to a gerbe on a Zariski open set $U\supset G^{reg}$
which meets each component of $Y$, so that its complement $Z$ has
codimension $\geq 2$. Also there will be a holomorphic connective structure on $\tcalC$.

The arguments which will lead to Theorem 2 are quite technical, as they
make heavy use of hypercohomology of a complex of sheaves with supports
in a closed subset, so many readers may wish to skip ahead to the statement of Theorem 2.

We denote by $V\subset G$ the set of elements $g$ of $G$ which
are not regular, or in other words $q^{-1}(g)$ is not finite. It is well-known that $V$ has
codimension $3$ in $G$. Then the obstruction to extending our gerbe with $0$-connection
from $G\min Z$ to $G\min V$ is an element of the hypercohomology group $H^3_{Z\min
V}(G\min V,\calO^*\to\Omega^1)$ with supports in $Z\min V$. Furthermore the
non-uniqueness of the extension is controlled by the hypercohomology group
$H^2_{Z\min V}(G\min V,\calO^*\to\Omega^1)$.
We will only say here a few
things about  hypercohomology  with supports, referring the reader to the book [{\bf K-S}]
for details. For any complex of sheaves $F^{\bu}$ over  a space $X$, and for $Y$ a closed
subset of $X$, the hypercohomology groups $H^p_Y(X,F^{\bu})$ with supports in $Y$ sit in
an exact sequence
$$\cdots\to H^{p-1}(X,F^{\bu})\to H^p_Y(X,F^{\bu})\to H^p(X,F^{\bu})
\to H^p(X\min Y,F^{\bu})\to\cdots$$
They can be computed by the \v Cech method using  open coverings $(V_\a)$ of $X$ 
and $(U_\a)$ of $X\min Y$ such
that $U_\a\subseteq V_\a$ and all cohomology groups $H^q(V_{\a_1\cdots \a_l},F^p)$ and
$H^q(U_{\a_1\cdots
\a_l},F^p)$ vanish for $q>0$. Then we can construct the \v Cech double complexes
$C^{\bu}(\calV,F^{\bu})$ and $C^{\bu}(\calU,F^{\bu})$ and we have a natural
restriction mapping from the first double complex to the second, which allows
to construct a triple complex, whose total cohomology is the hypercohomology
with supports. For the complex of sheaves $\calO^*\to\Omega^1$, we may pick the $V_\a$ and
$U_\a$ to be small open discs.

 We have an exact sequence of complexes of sheaves
$$0\to \BZ \to [\calO\to\Omega^1]\to[\calO^*\to\Omega^1]\to 0$$

Thus for cohomology with supports we have the exact sequence

$$0\to C\to H^3_{Z\min V}(G\min V,\calO^*\to\Omega^1)\to H^4_{Z\min V}
(G\min V,\BZ)$$
where $C$ is a complex vector space, namely a hypercohomology group
with coefficients in $\calO\to\Omega^1$. Now the group
$H^4_{Z\min V} (G\min V,\BZ)$ is the free abelian group generated by the cohomology
classes of those components $Z_j$ of $Z$ which have codimension $2$ in $G$. Taking
the inverse images
$\tilde V$,
$\tilde Z$ of $V$ and $Z$ in $\tilde G$, we have a similar exact sequence
$$0\to \tilde C\to H^3_{\tilde Z\min\tilde V}(\tilde G\min
\tilde V,\calO^*\to\Omega^1)\to H^4_{\tilde Z\min \tilde V} (\tilde G\min\tilde
V,\BZ)$$
Now the map $C\to\tilde C$ is injective because the mapping $\tilde G\min \tilde V
\to G\min V$ is finite and $C$ is a vector space.  The map
$H^4_{Z\min V}
(G\min V,\BZ)\to H^4_{\tilde Z\min \tilde V} (\tilde G\min\tilde
V,\BZ)$ is injective because the inverse image of a component $Z_j$ of $Z$
is a union of codimension $2$ components of $\tilde Z$. It follows that the
pull-back map
$H^3_{Z\min V}(G\min V,\calO^*\to\Omega^1)\to
H^3_{\tilde Z\min\tilde V}(\tilde G\min
\tilde V,\calO^*\to\Omega^1)$ is injective. Thus the obstruction to extending
our gerbe $\calC$ from $G\min Z$ to $G\min V$ vanishes, because
the gerbe $\tcalC$ over $\tilde G\min\tilde Z$ extends to $\tilde G\min
\tilde V$.

One shows similarly
 that the map
$H^2_{Z\min V}(G\min V,\calO^*\to\Omega^1)\to
H^2_{\tilde Z\min\tilde V}(\tilde G\min
\tilde V,\calO^*\to\Omega^1)$ is injective (this is actually easier, as these
cohomology groups are complex vector spaces). Thus the data of the gerbe
$\tcalC$ over $\tilde G\min
\tilde V$ leaves no amount of freedom for the extension of $\calC$
to $G\min V$.

Then using Lemma 6 we obtain

\theo{2}{Let $H$ be a connected reductive algebraic complex group,
and let $f:H\to G$ be an algebraic homomorphism where $H$ is also
algebraic reductive,  $Ker(f)$ is central, and $G=Z(G)\cdot Im(f)$.
For any $W$-invariant tensor
$b\in X^*(S)\ot X^*(T)$ such that the bilinear form $b:X_*(S)\ot
X_*(T)\to\BZ$ satisfies the condition (EV), the corresponding $H\times
W$-equivariant holomorphic gerbe $\tcalC$ with
$0$-connection over $\tilde G$ can be descended in an (essentially) unique way
to an $H$-equivariant gerbe $\calC$ over $G$.}

\cor{1}{If $inv:G\to G$ denotes the inverse map, then the pull-back
gerbe $inv^*\calC$ is equivalent to $\calC^{\ot-1}$.}

\proof{This follows easily from Lemma 3 and the fact that $\calC$ is obtained
from $\tcalC$ by quotienting by $W$ and then extending from an open set. \qed}

We denote by $\calC\boxtimes\calC$ (external tensor product) the gerbe over $G^2$ obtained by
tensoring the gerbes $p_1^*\calC$ and $p_2^*\calC$.
The following result is analogous to Proposition 3.2 in [{\bf A-M-M}].

\cor{2}{Assume $H=G$ so that $\calC$ is $G$-equivariant.
Let $D:G^2\to G^2$ be the double map
$D(a,b)=(ab,a^{-1}b^{-1})$. Then the pull-back of $\calC\boxtimes\calC$ under $D$ is a trivial
gerbe.}

\proof{The equivariance of $\calC$ implies that for $d_j:G^2\to G$ the face maps
of the Appendix, the tensor product gerbe $d_0^*\calC\ot d_1^*\calC^{\ot-1}$ is trivial.
Now using Lemma 3 this means that the gerbe $d_0^*\calC\ot d_1^*inv^*\calC$ is trivial.
If we introduce the mapping $\delta:G^2\to G^2$ such that
$\delta(a,b)=(d_0(a,b),d_1(a,b)^{-1})=(aba^{-1},b^{-1})$, we see that the pull-back gerbe
$\delta^*(\calC\boxtimes\calC)$ is trivial. Now we can write
$D=\delta \phi$, where $\phi(a,b)=(ab^{-1},ba)$, hence $D^*(\calC\boxtimes\calC)$ is trivial
too. \qed}

Here is a description of $\calC$: for an open subset $U$ of $G$, an object of $\calC_U$
is an object $P$ of $\tcalC_{f^{-1}(U)}$ together with the structure of a $W$-equivariant
object on the restriction $Q$ of $P$ to $\tcalC_{f^{-1}(U)\cap \tilde G^{reg}}$. This
means (see Appendix) that for any $w\in W$ we are given an isomorphism
$\eta_w:Q\tilde{\to}w^*(Q)$ such that

(1) the $\eta_w$ satisfy the cocycle condition
$\eta_{w_1w_2}=w_2^*(\eta_1)\eta_2$

(2) each $\eta_w$ has no poles along components of $f^{-1}(U)\cap (\tilde G\min \tilde
G^{reg})$. This means that if $w^k=1$, then $(w^{k-1})^*(\eta_w)(w^{k-2})^*(\eta_w)\cdots
\eta_w$ is an automorphism of $Q$ which extends holomorphically to an automorphism
of $P$ over $U$.

 Isomorphisms are isomorphisms of objects of
$\tcalC_{f^{-1}(U)}$ which are compatible with the  extra data $(\eta_w)$.
From this description of the gerbe $\tcalC$ one can show easily that it is $H$-equivariant.
The $0$-connection on $\calC$ can then be described in terms of an $\O^1$-torsor
$Co(P,\eta_w)$ attached to the data $(P,\eta_w)$: the sections of this sheaf are the
holomorphic sections $\na$ of $Co(P)$ over $U$ which are $W$-invariant in the sense that
 $\eta_w$ maps $\na$ to $w^*\na$.
 It would of
course be very nice to find some explicit gerbe data for the gerbe $\calC$.

\sec{6. Discussion of the combinatorial data.}

The bilinear forms $b:X_*(S)\ot X_*(T)$ satisfying the conditions of W-invariance
and (EV) were introduced independently by Toledano in [{\bf To}] and by the author and Deligne
in [{\bf B-D}].

For $G$ simply-connected and simple, and for $H=G$, the allowable $b\in X^*(T)\ot
X^*(T)$ are the integer multiples of the basic  $b_0$ for which $b_0(\check
\a,\check\a)=2$ for a long root
$\a$ (so that
$\check\a$ is a short coroot). This bilinear form $b_0$ is introduced in [{\bf P-S}]
for the purpose of constructing central extension of loop groups. 

For $G=SL(n,\BC)$, $T$ the group of diagonal matrices, $X^*(T)$
is the quotient of the free group on the diagonal entries
$(t_1,\cdots,t_n)$ by the relation $t_1+\cdots+t_n=0$.
The roots are the $\a_{ij}=t_i-t_j$ for $i\neq j$.
Let $(e_1,\cdots,e_n)$ be the basis dual to $(t_1,\cdots,t_n)$.
The dual group $X_*(T)$ is the subgroup of $\BZ e_1\cdots\BZ e_n$
comprised of the linear combinations $\sum n_ie_i$ such that $\sum n_i=0$.
The coroots are the $\check \a_{ij}=e_i-e_j$. The basic element
$b_0$ is $b_0=\sum_{i=1}^n~ t_i^2$.

Now take $G_{ad}=PGL(n,\BC)$ to be the adjoint group of $SL(n,\BC)$, so that $T$ is
replaced by the quotient group
$T_{ad}=(\BC^*)^n/\BC^*_{\rm diag}$, where $\BC^*$ is embedded diagonally in
$(\BC^*)^n$. Then $X^*(T_{ad})=Q$ is the coroot lattice generated by the
$\check\a_{ij}$. If we take $H=G_{ad}=PGL(n,\BC)$ then the allowable $b\in
X^*(T_{ad})\ot X^*(T_{ad})$ are the integer multiples of
$\sum_{i<j}~\check\a_{ij}^2=n\cdot b_0$. The same situation occurs if we take $H=G$
mapping to $G_{ad}$ in the obvious way. Now let us consider $H=G=GL(n,\BC)$, and
take $T$ to be the group all diagonal matrices, so that
$X^*(T)=\BZ t_1\oplus\cdots\oplus\BZ t_n$. 
If we take $H=G$, then the allowable $b$ are of the form
$b=l\sum_{i=1}^n t_i^2+m (t_1+\cdots+t_n)^2$, for $l,m\in\BZ$.
If now $H=SL(n,\BC)$, so that $X^*(S)$ is the quotient of $X^*(T)$ by the relation
$t_1+\cdots+t_n=0$, then we are left with the integer multiples of $t_1^2+\cdots+t_n^2$.

If we take $G$ to be of type $B_n$, then for the simply-connected group
$G_{sc}=Spin(2n+1,\BC)$ the character group $X^*(T_{sc})$ is the group
$1/2\sum_i m_it_i$ where $m_i\in\BZ$ and $m_1\equiv\cdots\equiv m_n \bmod 2$. The simple
coroots are $(e_1-e_2,\cdots,e_{n-1}-e_n,2e_n)$.
For $H=G=Spin(2n+1,\BC)$, the
allowable
$b$ are the integer multiples of $b_0=\sum_{i=1}^nt_i^2$.
For the adjoint group $G_{ad}=SO(2n+1,\BC)$ and for $H=G_{ad}$ or $H=G_{sc}$,
the allowable $b$ are integer multiples of $2b_0$. This corresponds to the
well-known phenomenon that only the tensor square of the fundamental gerbe
over $Spin(2n+1)$ can be descended to $SO(2n+1)$.

If we take $G$ to be of type $D_n$, then for $G_{sc}=Spin(2n,\BC)$
the character group $X^*(T_{sc})$ is the group of $1/2\sum_{i=1}^n m_it_i$,where
$m_1\equiv m_2\cdots\equiv m_n\bmod 2$. The coroots are the $e_i-e_j$.
Thus for $H=G=G_{sc}$, the allowable $b$ are the integer multiples of
$b_0=\sum_{i=1}^n t_i^2$. Now take $G=SO(2n,\BC)$ so that $X^*(T)=\sum_{i=1}^n
\BZ t_i$; then for $H=G_{sc}$ or $H=Spin(2n,\BC)$ we find the allowable $b$ are
again the integer multiples of $b_0$. Now if we take $G=G_{ad}=SO(2n,\BC)/\{ \pm
1\}$, then
$X^*(T_{ad})$ is the group of $\sum m_it_i$ such that $m_1+\cdots+ m_n$ is even.
Then for any choice of $H$ (equal to $G_{sc}$ divided by a central subgroup),
we find that only integer multiples of $2b_0$ is allowable.

\sec{7. Restricting the gerbe to conjugation orbits.}

\lem{7}{For any regular semisimple element $g$ of $G$, the restriction
of the $H$-equivariant gerbe $\calC$ to the $H$-orbit $\calO_g=H\cdot g\subset G$
is $H$-equivariantly trivial.}

\proof{First of all we can lift $\calO$ to an  $H$-orbit $\tilde\calO$ in
$\tilde G$ which maps isomorphically to $\calO$, namely the orbit of $(g,B)$ for any Borel
subgroup of $G$ containing
$g$. Then it is enough to show that the restriction of the gerbe $\tcalC$ to $\tilde\calO$
is $H$-equivariantly trivial. Using the description of the gerbe $\tcalC$
given in \S 3, we see that any logarithm $\xi$ of $g$ gives such an equivariant
trivialization. \qed}

In general, the restriction of $\calC$ to an $H$-orbit will not be equivariantly
trivial. As shown in the Appendix, the obstruction to the triviality is a central
extension of the centralizer group. We will describe this central extension
in case $g$ is semisimple. First we recall some well-known facts about central extensions
of reductive complex groups. 
Let $L$ be a reductive connected algebraic group over
$\BC$. We want  to describe in combinatorial terms the
 central extensions $1\to \BC^*\to \tilde L\to L\to 1$ of complex
algebraic groups. Let $S$ be a maximal torus in $L$. The inverse image
$\tilde S$ is a maximal torus in $\tilde L$, and we have an extension
of free abelian groups
$$0\to X^*(S)\to X^*(\tilde S)\to X^*(\BC^*)=\BZ\to 0.$$
Hence the extension (EXT) is equivariant under the Weyl group of $S$ in $L$, which we will
simply denote by $W_L$.

\lem{8}{The central extension $\tilde L$ of $L$ is entirely determined by
the extension (EXT) of $W_L$-modules.}

Now the extension (EXT) gives rise to a degree $1$-cocycle group cocycle
$w\mapsto c_w:W_L\to X^*(S)$ as follows: pick an element $\chi\in X^*(\tilde S)$ which
restricts to $1\in X^*(\BC^*)$. Then $c_w=w\chi-\chi\in X^*(S)$ is a $1$-cocycle, and its
cohomology class vanishes iff the extension (EXT) has a $W_L$-equivariant
splitting.

Let us apply this to  $f:H\to G$ as in \S 2 and some semisimple
element $g\in G$. 
Denote by $L$ the connected component of the centralizer
of $g$ in $H$. Then pick a maximal torus $T$ of $G$ containing $g$ and as usual let
$S=f^{-1}(T)\subset H$, which is a maximal torus of $L$. The Weyl group
$W_L$ of $S$ in $L$ is the subgroup of $W=W_H$ which centralizes
$g$. Pick some   element $\xi$ of $\ft$ such that $exp(2\pi i\xi)=g$. Then
for $w\in W_L$, the difference $w\xi-\xi$ has exponential $1$,
 so $d_w=w\xi-\xi$ belongs to $X_*(T)$. This gives a $1$-cocycle of
$W_{S,L}$ with values in $X_*(T)$.

Now given $b\in X^*(S)\otimes X^*(T)$ which is $W$-invariant
and satisfies (EV), we have constructed an $H$-equivariant holomorphic gerbe
over $G$. Recall that $\un b:X_*(T)\to X^*(S)$ is the corresponding linear map.
Then according to the Appendix, for any $g$ we have a central extension
of the centralizer group $Z_g(H)$. We can now identify this central
extension:

\prop{6}{For $g\in G$ semisimple, the central extension of the centralizer
$L$ of $g\in H$ is described by the cohomology class of the group $1$-cocycle
$$w\in W_L\mapsto \un b(d_w)\in X^*(S).$$}

\proof{It follows from the Appendix that
the central extension $\tilde L$ of $L$ has restriction to $S$ given by the multiplicative
$\BC^*$-bundle $\ot_{i=1}^r (\chi_i,\ze_i(g))$ if $b\in X^*(S)\ot X^*(T)$ is equal
to $\sum \chi_i\ot\ze_i$. Here $\ze_i(g)$ is just a complex number and $\chi_i$ is a
holomorphic function $S\to\BC^*$. The expression $(\chi_i,\ze_i(g))$ denotes a Deligne
line bundle. A multiplicative trivialization of the Deligne line bundle above is obtained
using the logarithm $d\ze_i(\xi)$ of $\ze_i(g)$. This gives  a trivialization
$s=\ot_i (\chi,d\ze_i(\xi)]$ of the central extension of $S$, which however is not
$W_L$-equivariant. For
$w\in W_L$, we can measure the defect of equivariance as a group homomorphism
$c_w: S\to\BC^*$ given by $c_w(h)=ws(w^{-1}h)-s(h)$. Now we have
$$w\cdot s=\ot_i(w\chi_i,d\ze_i(\xi))=\ot_i(\chi_i,dw^{-1}\ze_i(\xi)]
=\ot_i(\chi_i,d\ze_i(w\xi)]$$
using the $W$-invariance of $b$.
It follows that
$$c_w(h)=\prod_i \chi_i(h)^{d\ze_i(d_w)}=[\un b(d_w)](h).$$ \qed}

It appears likely that this central extension is trivial for $SL(n,\BC)$
or $GL(n,\BC)$ but often non-trivial for other groups. In any case, its
cohomology class has finite order (bounded by the order of $W_L$), so if we replace the gerbe
by some tensor power the obstruction will vanish. This fits with the result
in [{\bf G-H-J-W}] and [{\bf A-M-M}] that the class in $H^3_K(K,\BR)$ restricts
trivially to any orbit.

We give an example where the central extension is non-trivial.

Pick for instance $H=G=Spin(7,\BC)$ and $g=exp(2\pi i\xi)$ where $\xi={e_1-e_2\over 2}$.
The centralizer $L$ of $g$ is infinitesimally isomorphic to $SL(2,\BC)^3$.
Then the subgroup $W_L$ of the Weyl group is the group $(\BZ/2)^3$
generated by the change of sign $(-1,-1,1)$, the change of sign $(1,1,-1)$ and by the
permutation of $1$ and $2$. The $1$-cocycle $d_w$ of $W_L$ with values in $X_*(T)$ is
$d_w=w\xi-\xi$. We now take $b=b_0=t_1^2+t_2^2+t_3^2$; then $\un b:X_*(T)\to X^*(T)$
is the inclusion which maps $e_i$ to $t_i$. Now we claim that $\un b(d_w)$ is not the
coboundary of some $u\in X^*(T)$. Since there is no vector in $X^*(T)$ fixed by
the group $W_L$, the only possible choice of $u$ would be $u=\un b(\xi)={t_1-t_2\over 2}$.
However this is a vector in
$X^*(T)\ot\BQ$ which does not belong to $X^*(T)$. Thus $\un b(d_w)$ is not a
coboundary and the corresponding central extension $\tilde L$ of $L$ is not trivial.
Clearly by construction it has order $2$, hence is induced by a central extension of
$L$ by $\BZ/2$, i.e., a connected double cover $\hat L$ of $L$. We can describe $L$ more
precisely as follows: $W_L$ is generated by the reflections corresponding to the
orthogonal roots $t_1-t_2,t_1+t_2,t_3$. For the maximal torus $S$ in $(SL(2,\BC))^3$,
$X_*(S)$ is spanned by the corresponding coroots $e_1-e_2,e_1+e_2,2e_3$. This subgroup is
on index $2$ in $X_*(T)$, and the quotient subgroup is spanned by the half-sum
$1/2((e_1-e_2)+(e_1+e_2)+2e_3)=e_1+e_3$. It follows that $L$ is isomorphic to the quotient
of $(SL(2,\BC))^3$ by the group $\BZ/2$ embedded diagonally in the center $(\BZ/2)^3$, and
thus $\hat L$ is the universal cover $(SL(2,\BC))^3$.

\sec{Appendix. Equivariant gerbes}

1) \und{Equivariant line bundles}

We will review well-known material on equivariant line bundles in a framework which will
adapt to gerbes. Let $G$ be a Lie group which acts smoothly on the smooth manifold $M$.
Let $m:G\times M\to M$ denote the action and let
$p_2:G\times M\to M$ be the second projection.
First of all if we view a line bundle as given by transition
cocycles $g_{\a\b}$ with respect to an open covering $V_\a$, then the problem we
encounter if that we  can't in general assume that each $V_\a$ is $G$-invariant.
In order to write down the equivariance data for the line bundle in \v Cech form, what we
need to do is cover
$G\times M$ by the open sets
$Z_{\a\b}=m^{-1}(G\times V_\a)\cap p^{-1}(G\times V_\b)$. Then to make our line bundle
$G$-equivariant we need to introduce a function
$\phi_{\a\b}:Z_{\a\b}\to\BC^*$. The meaning of $\phi_{\a\b}$ is as follows: if
$s_\a$ is the chosen non-vanishing section of $L$
over $V_\a$, then if $L$ is $G$-equivariant we have the section $(g,x)\mapsto
[g^*s_\a](x)$ over $m^{-1}(G\times
V_\a)$, and $\phi_{\a\b}$ is this section divided by $s_\b(x)$. Then we see that
$\phi_{\a\b}$ must satisfy the following requirement:

$$\phi_{\g\b}=[m^*g_{\g\a}]\phi_{\a\b}~~,~~\phi_{\a\delta}=[p_2^*g_{\delta\b}^{-1}]\phi_{\a\b}.$$
over the relevant open sets in $G\times M$.

Next a connection on $L$ amounts to $1$-forms $A_\a$ over $V_\a$ such that
$A_\b-A_\a=d\log (g_{\a\b})$. The connection is $G$-equivariant iff
we have the equality of $1$-forms on $G\times M$ in the direction of $M$:
$$m^*A_\a-p_2^*A_\b=d\log\phi_{\a\b}.$$

Of course, rather than using such data, it may be more convenient to use the geometric
language of line bundles and their pull-backs.
For this purpose, we introduce the simplicial manifold $G^{\bu}\times M$, which
is the family of manifolds
$G^n\times M$, equipped with the following face maps $d_i:G^n\times M\to G^{n-1}\times M$:

$$\matrix{d_0(g_1,\cdots,g_n,x)&=&(g_2,\cdots,g_n,g_1x)\hfill\cr
d_i(g_1,\cdots,g_n,x)&=&(g_1,\cdots,g_{i-1},g_ig_{i+1},\cdots,g_n,x)~
\rm{for}~1\leq i\leq n-1\cr
d_n(g_1,\cdots,g_n,x)&=&(g_1,\cdots,g_{n-1},x)\hfill}$$

The significance of the simplicial manifold $G^{\bu}\times M$ is that
it is a geometric model for the Borel space $EG\times^GM$. Recall that the
cohomology of the Borel space is the (Borel) equivariant cohomology
$H^*_G(M)$.

Note that $d_0:G\times M$ is equal to $m$, and $d_1:G\times M\to M$ is equal to
$p_2$.

Then an equivariant line bundle over $M$ is a line bundle $L$ over $M$, equipped with
a non-vanishing section $\s$ of the line bundle $d_0^*L\ot d_1^*L^{\ot-1}$,
which satisfies the cocycle condition
$d_0^*\s\ot d_1^*\s^{\ot-1}\ot d_2^*\s=1$.
This makes sense as $d_0^*\s\ot d_1^*\s^{\ot-1}\ot d_2^*\s$
is a section of the trivial line bundle over $G^2\times M$ (due to the relations among
the  iterated face maps $G^2\times M\to M$).

When $G$ is discrete, $\s$ amounts to a family of isomorphisms
$\s_g:L\tilde{\to} g^*L$, and the cocycle condition becomes simply
$\s_{g_1g_2}=g_2^*\s_{g_1}\s_{g_2}$. For a Lie group $G$, $\s$ still amounts to the data
of such $\s_g$, which must vary smoothly as a function of $g\in G$.

From such $\s$ we obtain the previous function 
$\display\phi_{\a\b}={[d_0^*s_\a\ot p_2^*s_\b^{\ot-1}]\over\s}$
- a ratio of two sections of $d_0^*L\ot d_1^*L^{\ot-1}$ over $Z_{\a\b}$.

Then the condition that a connection $\na$ is $G$-invariant amounts to the
condition that the covariant derivative of $\s$ vanishes in the $M$-direction.

Recall  that line bundles over $M$ are classified by the cohomology group
$H^2(M,\BZ)$. Similarly we have

\prop{A-1}{For $G$ a compact Lie group, the $G$-equivariant line bundles over $M$ 
are classified by the (Borel) equivariant cohomology group $H^2_G(M,\BZ)$.}

\proof{The data $(g_{\a\b},\phi_{\a\b})$ can viewed as a \v Cech cocycle for
the simplicial manifold $G^n\times M$, when we use the covering $(V_\a)$ of 
$M$, the open covering $Z_{\a\b}$ of $G\times M$ and the open covering
$(d_0^*V_\a)\cap (d_1^*V_\b)\cap (d_2^*V_\g)$ of $G^2\times M$. Thus it yields a
class in the degree $1$ hypercohomology $H^1(G^{\bu}\times M,\un{\BC}^*)$. Now the
exponential exact sequence gives an exact sequence
$$H^1(G^{\bu}\times M,\un{\BC})\to H^1(G^{\bu}\times M,\un{\BC}^*)
\to H^2(G^{\bu}\times M,\un{\BZ})\to H^2(G^{\bu}\times M,\un{\BC})$$
Now all the  hypercohomology groups
$H^p(G^{\bu}\times M,\un{\BC})$ can be computed in terms of the global sections
of the sheaves $\un\BC$ over $G^k\times M$, since $\un\BC$ is a fine sheaf.
Then the complex in question becomes the complex of smooth cochains of $G$ with
values in the $G$-module $C^{\inf}(M)$. Thus the hypercohomology groups 
$H^p(G^{\bu}\times M,\un{\BC})$ are exactly the differentiable cohomology
groups $H^p(G,C^{\inf}(M))$, which are $0$ for $p>0$ since $G$ is compact.
Thus $H^1(G^{\bu}\times M,\un{\BC}^*)$ identifies with 
$H^2(G^{\bu}\times M,\un{\BZ})$, and the latter group is just the \v Cech version 
of $H^2_G(M,\BZ)$. \qed}

It is interesting to see more concretely how the data
$(g_{\a\b},\phi_{\a\b},A_\a)$ for an equivariant line bundle lead to an equivariantly
closed differential form $(F,\mu)$. First of all $F$ is the curvature so that
$F_{/V_\a}=dA_\a$. Next we consider the $1$-form
$\o_{\a\b}$ on $Z_{\a\b}\subset G\times M$ defined as
$$\o_{\a\b}=d_0^*A_\a-d_1^*A_\b-d \log \phi_{\a\b}.$$
This $1$-form vanishes in the $M$-direction, so is entirely
in the $G$-direction. Also $\o_{\a\b}$ is $G$-invariant if $G$ acts on $G$ by left
multiplication. Further over $Z_{\a\b}\cap Z_{\g\b}$ we have
 $$\align{\o_{\g\b}-\o_{\a\b}&=m^*(A_\g-A_\a)-
d \log({\phi_{\g\b}\over \phi_{\a\b}})\cr
&=d\log (m^*g_{\g\a})-
d \log({\phi_{\g\b}\over \phi_{\a\b}})=0}$$
using (A-1) and similarly $\o_{\a\b}$ coincides with $\o_{\a\delta}$ over
$Z_{\a\b}\cap Z_{\a\delta}$. Hence the $\o_{\a\b}$ glue together to give  a global
$1$-form
$\o$ on $G\times M$ which is $G$-invariant and lives in the $G$-direction. We can 
write
$\o=\lan p_1^*g^{-1}dg,\mu\ran$, where $g^{-1}dg$ is the Maurer-Cartan $1$-form on
$G$ and $\mu:M\to\fg^*$ is a smooth function. This function is a moment map
for the $G$-action. We can evaluate it as follows: for $\xi\in\fg$, denote by
$\tilde\xi$ the corresponding vector field on $M$, and by $(\xi,0)$ the 
vector field on $G\times M$ which lives in the $G$-direction and is the 
left-invariant vector field defined by $\xi$. The derivative
$[(\xi,0)]\cdot \log \phi_{\a\b}]_{(1,x)}$ is equal to $\display\xi\cdot s_\a\over
s_\a$, where $\xi\cdot s_\a$ denotes the derivative at $t=0$ of
$exp(t\xi)^*\cdot s_\a$. Then we find:
$$\lan \mu(x),\xi\ran=\lan \o,[(\xi,0)]\ran=
{\na_{\tilde\xi}s_\a\over s_\a}-{\xi\cdot s_\a\over s_\a}$$
which is a standard description of the moment map as measuring the difference 
between two infinitesimal actions of $\fg$ on sections of $L$: the one given by the
connection evaluated along the $G$-orbits
and the one given by the $G$-action on sections of the equivariant line bundle $L$
[{\bf B-V}].
One checks easily that $d\lan\mu,\xi\ran=\lan\tilde\xi, F\ran$ as required, so that
$(F,\mu)$ is an equivariantly closed $2$-form. This is the equivariant Chern class 
as constructed by Berline and Vergne [{\bf B-V}].

2) \und{Equivariant gerbes}

We will discuss equivariant gerbes in  a similar spirit as we discussed equivariant
line bundles. Let
$(V_\a,\La_{\a\b},\th_{\a\b\g})$ be some gerbe data over $M$ (as mentioned before, 
we will suppress from the notation the isomorphism between $\La_{\a\b}^{\ot-1}$ and
$\La_{\b\a}$). Then to make the gerbe data
$G$-equivariant we need to pick a line bundle $E_{\a\b}$ over 
$Z_{\a\b}\subseteq G\times
M$ together with isomorphisms

$$\phi_{{\g/\a},\b}:d_1^*\La_{\g\a}\ot E_{\a\b}\tilde{\to}E_{\g\b}$$
and
$$\phi_{\a,{\delta/\b}}:
d_0^*\La_{\delta\b}^{\ot-1}\ot E_{\a\b}\tilde{\to}E_{\a\delta},$$
which satisfy the compatibility conditions
$$\phi_{{\epsilon/\a},\b}=d_1^*\th_{\epsilon\g\a}\ot[\phi_{{\epsilon/\g},\b}
\phi_{{\g/\a},\b}],$$
$$\phi_{\a,{\epsilon/\b}}=d_0^*\th_{\epsilon\delta\b}^{\ot-1}
[\phi_{\a,{\epsilon/\delta}}
\phi_{\a,{\delta/\b}}]$$
and the obvious commutation relation between the two types of isomorphisms.

Next we need a non-vanishing section $\psi_{\a\b\g}$
of the line bundle
$Q_{\a\b\g}=d_0^*E_{\b\g}\ot d_1^*E_{\a\g}^{\ot-1}
\ot d_2^*E_{\a\b}$ over $d_0^*Z_{\b\g}\cap d_1^* Z_{\a\g}
\cap d_2^* Z_{\a\b}\subset G^2\times M$. This section should satisfy three 
conditions: first of all, $\psi_{\a\b\g}$ should correspond to $\psi_{\delta\b\g}$
under the tensor product of the isomorphisms $d_1^*\phi_{\delta/\a,\g}$
and $d_2^*\phi_{\delta/\a,\b}^{\ot-1}$; there are two similar conditions involving
changing the second and third indices in $Q_{\a\b\g}$.
Secondly, we require the cocycle condition
$$d_0^*\psi_{\b\g\delta}\ot d_1^*\psi_{\a\g\delta}^{\ot-1}\ot
d_2^*\psi_{\a\b\delta}
\ot d_3^*\psi_{\a\b\g}^{\ot-1}=1.$$
This makes sense as the left hand side is a section of the trivial line bundle
over an open set of $G^3\times M$.

We can now state

\prop{A-2}{If $G$ is a compact Lie group, the equivariant DD-gerbes over $M$ are
classified by the equivariant cohomology group $H^3_G(M,\BZ)$.}
\proof{
The proof is similar to that of Proposition A-1 so we will be brief. One first assumes that all
line bundles
$\La_{\a\b}$ and $E_{a\b}$ are trivial. Then one can interpret the data
$(\th_{a\b\g},\phi{\g/\a,\b}\phi_{\a,\delta/\b},\psi_{\a\b\g})$ as yielding as giving
a \v Cech $2$-cocycle with values in the simplicial sheaf $\un{\BC}^*$ on the
simplicial manifold $G^{\bu}\times M$. Then one uses the exponential exact sequence
to compare $H^2(G^{\bu}\times M,\un\BC^*)$ with $H^3(G^{\bu}\times M,\un\BZ)$. \qed}

There is one case where these unwieldy constructions are simpler: assume that all open sets
$V_\a$ are
$G$-stable, so that $d_0^{-1}(V_\a)=d_0^{-1}(V_\a)$ and that the $\La_{\a\b}$ are
equivariant line bundles over $V_{\a\b}$ and that each $\th_{\a\b\g}$ is $G$-invariant.
Then we can pick $E_{\a\b}=d_0^*\La_{\a\b}$; then we can take $\phi_{\a,\delta/\b}$
to be induced by $d_0^*\th{\b\a\delta}$, and since $d_1^*\La_{\a\b}$ is isomorphic
to $d_0^*\La_{\a\b}$ (by the equivariance of $\La_{\a\b}$), we can take
$\phi_{\g/\a,\b}$ to be the isomorphism induced by $d_1^*\th_{\g\b\a}$. Then
$\psi_{\a\b\g}$ is induced by $\th{\a\b\g}$ and all conditions are satisfied
for an equivariant gerbe.

 The data also simplify considerably in
the case where $M$ is a point. Then there is no need for the covering
$(V_\a)$ and the line bundle $\La_{\a\b}$, so the construction boils down to
a line bundle $E$ over $G$ and a non-vanishing section $\psi$ of the line bundle
$d_0^*E\ot d_1^*E^{\ot-1}\ot d_2^*E$ over $G^2$. The fiber of this line
bundle at $(g_1,g_2)$ is equal to $E_{g_2}\ot E_{g_1g_2}^{\ot -1}
\ot E_{g_1}$. Thus $\psi$ amounts to a multiplicative structure on the line
bundle
$E\to G$, namely an isomorphism $E_{g_1}\ot E_{g_2}\tilde{\to} E_{g_1g_2}$
which varies smoothly with $(g_1,g_2)$. This means also that the total space 
$\tilde G$ of the corresponding $\BC^*$-bundle acquires a product structure
$\tilde G\times \tilde G\to\tilde G$ which lifts the product on $G$.
The cocycle condition for $\psi$ just says that this product law is associative.
Then $\tilde G$ becomes a Lie group, which is a central extension
of $G$ by $\BC^*$. In fact we have

\prop{A-3}{The equivalence classes of $G$-equivariant gerbes over a homogeneous 
space
$G/H$ correspond precisely to the central extensions
$$1\to\BC^*\to\tilde H\to H\to 1.$$}

We have proved this in the case where $G=H$. We next claim that by restriction
to the base point of $G/H$, we get a bijective correspondence between 
$G$-equivariant gerbes over $G/H$ and $H$-equivariant gerbes over a point. What we
need to do is to give an analog for gerbes of the construction of a homogeneous
vector bundle over a homogeneous space. So we start with an $H$-homogeneous gerbe
$\calC$  over a point, and we pull it back to a gerbe over $G$ which
is left
$G$-equivariant and right $H$-equivariant. Then we have the following lemma:

\lem{A-1}{If the Lie group $H$ acts freely on a manifold $X$, then pull-back gives 
a  bijective correspondence between equivalence of gerbes on $X/H$ and equivalence
classes of $H$-equivariant gerbes over $X$.}

Thus we can descend our gerbe to a gerbe on $G/H$, which is still $G$-equivariant
as the left action of $G$ on $G$ commutes with the right $H$-action.
This ``homogeneous gerbe construction'' is inverse to restriction to the base point.
\qed

This construction can be generalized to the case where the gerbe data is trivial, 
so that the line bundles $\La_{\a\b}$ are trivial and $\th_{\a\b\g}=1$. Then the
isomorphisms $\phi_{\g/\a,\b}$ give, for fixed $\b$, gluing isomorphisms between
the line bundles $E_{\a\b}$ and $E_{\g\b}$ over the overlap $Z_{\a\b}\cap Z_{\g\b}$
of their domains of definition. Then we can use these gluing data to produce a line
bundle $E_\b$ over $d_1^*V_\b$. The isomorphisms $\phi_{\a,\delta/\b}$ then glue
together, as $\a$ varies, to yield
a global isomorphism between $E_\b$ and $E_\delta$ (this uses the commutation
relations between the two types of isomorphisms). So finally we have a global line
bundle
$E$ over $ G\times M$. Then the $\psi_{\a\b\g}$ glue together to give a global
non-vanishing section $\psi$ of $d_0^*E\ot d_1^*E^{\ot-1}\ot d_2^*E$
over $G^2\times M$. To interpret this, recall that $A=G\times M$ has the structure
of a differentiable groupoid, where $(g,x)$ is viewed as the arrow labeled by $g$
which goes from $x$ to $g\cdot x$. Then $\psi$ gives the total space $\tilde A$
of the $\BC^*$-bundle associated to $E$ a partial composition law which lifts that
on the groupoid $A$. The cocycle condition for $\psi$ means that this composition
law is associative, so that $\tilde A$ becomes a groupoid which is a central
extension of $A$ by $\BC^*$ (we refer to [{\bf W-X}] for central extensions
of groupoids and their applications to geometric quantization).

Next we examine the notion of equivariant $0$-connection on an equivariant gerbe.
Let then $(D_{\a\b})$ be a $0$-connection. To make it $G$-equivariant, we need to
pick a connection $\na_{\a\b}$  in the direction of $M$ on each $E_{\a\b}$; this
means that the covariant derivative of a section of $E_{\a\b}$ on an open set of
$G\times M$ is a $1$-form in the $M$-direction. Such a connection is also called a
relative connection (relative to the projection $G\times M\to M$). We require that
the isomorphisms
$\phi_{\g/\a,\b}$ and $\phi_{\a,\delta/\b}$ are compatible with the relative
connections. Then $\psi_{\a\b\g}$ should be a horizontal section of $Q_{\a\b\g}$,
but only in the direction of $M$.

Given a $1$-connection $(F_\a)$, it is equivariant iff it satisfies the constraint
 $$Curv(\na_{\a\b})=d_0^*F_\a-d_1^*F_\b~
{\rm ~in~the~direction~of~}M~ {\rm (over}
~Z_{\a\b})$$ Then the $3$-curvature $\O$ satisfies $d_0^*\O=d_1^*\O$, so it is
$G$-invariant. It is most natural at this point to write down the equivariantly 
closed
$3$-form which is the equivariant Chern character of the equivariant gerbe. As in 
the case of line bundles, we introduce the $2$-form
$d_0^*F_\a-d_1^*F_\b-Curv(\na_{\a\b})$ on $Z_{\a\b}$ which has zero component in
$A^0(G)\hat\ot A^2(M)$ . We look at the component 
$B_{\a\b}$ of this $2$-form
onto the factor $A^1(G)\hat\ot A^1(M)$ of $A^2(G\times M)$. It is easy to see that
the $B_{\a\b}$ glue together to give a global $2$-form $B$ on $G\times M$.
As $B$ is $G$-invariant
for the left action of $G$ on itself, it can be written down
as $B=\lan g^{-1}dg,\ot E\ran$, where $E$ is a $\fg$-valued $1$-form on $M$
and $B$ is obtained using the evaluation map 
$\lan\pha{a},\pha{a}\ran:\fg^*\ot\fg\to\BR$.
Then $(\O,E)$ is an equivariantly closed $3$-form on $M$.

Now we interpret the notion of equivariant gerbe data in terms of DD-gerbes.
So let $\calC$ be a DD-gerbe over $M$, viewed as in \S 1 as a sheaf of groupoids
satisfying axioms 1)-3). Then to make $\calC$ $G$-equivariant we need two extra
pieces of data, First we need an equivalence $\phi:d_1^*\calC\to d_0^*\calC$ of
gerbes over $G\times M$. Equivalently, $\phi$ amounts to a global object $R$ of the  gerbe
$d_0^*\calC\ot d_1^*\calC^{\ot-1}$ over $G\times M$. Second we need an isomorphism
$$\psi:d_0^*R\ot d_1^*R^{\ot -1} \ot d_2^*R  
\tilde{\to}\bf 1$$
of objects of the trivial gerbe over $G^2\times M$.
This isomorphism must satisfy the cocycle condition
$$d_0^*\psi\ot d_1^*\psi^{\ot-1}\ot d_2^*\psi\ot d_3^*\psi^{\ot-1}=1.$$

For $G$ discrete, $\phi$ amounts to gerbe equivalences $\phi_g:\calC\to g^*\calC$
and $\psi$ amounts to natural transformations 
$$\psi_{g_1,g_2}:\phi_{g_1g_2}\to (g_2^*\phi_{g_1})\phi_{g_2}$$
between equivalence of gerbes
which must satisfy a cocycle condition (that condition may be visualized as a commutative
tetrahedron). We can write $\phi_g(P)$ as $g_*P$ for an object $P$ of $\calC$
over some open set. Then $\psi_{g_1,g_2}$ is an isomorphism between
$(g_1g_2)_*P$ and $(g_1)_*[(g_2)_*P]$.

We briefly adumbrate how these data lead to the equivariant gerbe data
discussed previously. First take an open covering $(V_\a)$ of $M$ and objects
$P_\a\in\calC_{V_\a}$. Then as in \S 1, we have the line bundle $\La_{\a\b}
\to V_{\a\b}$ associated to the $\BC^*$-bundle $\un{Isom}(P_\b,P_\a)$.
Over $Z_{\a\b}$ we have the objects $\phi(d_1^*P_\b)$ and $d_0^*P_\a$ of the
pull-back gerbe $d_0^*\calC$. Then we define the line bundle $E_{\a\b}$ to be the
line bundle associated to the
$\BC^*$-bundle $\un{Isom}(d_0^*P_\b,\phi(d_1^*P_\a))$. The isomorphisms
$\phi_{\g/\a,\b}$ and $\phi_{\a,\delta/\b}$ are given by composition of isomorphisms
of objects in the gerbe $d_0^*\calC$. Then the tensor product line bundle
$Q_{\a\b\g}$ corresponds to the tensor product $\BC^*$-bundle given by 
$$\un{Isom}(pr_2^*P_\g,pr_1^*P_\b)\otimes\un{Isom}(pr_0^*P_\a,pr_2^*P_\g)
\otimes\un{Isom}(pr_1^*P_\b,pr_0^*P_\a).$$
Thus by composition of isomorphisms we obtain the trivialization
$\psi_{\a\b\g}$ of $Q_{\a\b\g}$, which by its construction satisfies a cocycle
condition.

Another advantage of the notion of equivariant gerbe $\calC$ is that we can define an
equivariant object of $\calC$. This means that $P$ is an object
of $\calC_M$, equipped with an isomorphism $\eta:\phi(d_1^*P)\tilde{\to}d_0^*P$
of objects of $[d_0^*\calC]_{G\times M}$ which satisfies the associativity condition
$d_0^*\eta\ot d_2^*\eta^{\ot-1}\ot d_2^*\eta=1$. For $G$ discrete, $\eta$ amounts to a
family of isomorphisms $\eta_g:\phi_g(P)\tilde{\to} g^*P$ which satisfy
 the cocycle condition. This is formally the same description as for
equivariant line bundles.

In case $M$ is a point, the obstruction to finding an equivariant object
of $\calC$ is a central extension of $G$ by $\BC^*$. Indeed, picking any
$P$ and $\eta$ as above, the automorphism $d_0^*\eta\ot d_2^*\eta^{\ot-1}\ot
d_2^*\eta$ of $pr_1^*P\in pr_1^*\calC_{G^2}$ is a function $G^2\to\BC^*$, which is
a $2$-cocycle. Thus we recover the central extension we previously described using
equivariant gerbe data.

Now given a $0$-connection on the gerbe $\calC$, viewed as in \S 1 as a sheaf
$Co(P)$ attached to each object of each $\calC_U$, to make the $0$-connection
equivariant we need to extend the equivalence $\phi:d_1^*\calC\to d_0^*\calC$
of gerbes over $G\times M$ to an equivalence of gerbes with $0$-connection.
This means that for each object $P$ of each $(d_1^*\calC)_U$ we give an isomorphism
of sheaves $\phi_*:Co_{rel}(P)\to Co_{rel}(\phi_* P)$, where $Co_{rel}(P)$ is the sheaf
obtained by dividing $Co(P)$ by the action of the $1$-forms on $G\times M$ which are
in the direction of $G$. 
Thus $Co_{rel}(P)$ is a torsor under the sheaf $\O^1_{M\times G\to G}$ of relative $1$-forms
with respect to the projection $G\times M\to G$. 
In case $G$ is  discrete,
we have $\O^1_{M\times G\to G}=\O^1_{M\times G}$, and a section of this sheaf is a family
$\o_g$ of $1$-forms on (open sets of) $M$, indexed by $g\in G$. Let us see what $\phi_*$ looks
like when $P=d_1^*Q$ where $Q$ is some object of $\calC$ over some open set. Then
we put $\phi_g(Q)=g_*(Q)$ as explained earlier, so that $\phi_*$ amounts to
a family of isomorphisms of torsors $Co(Q)\tilde{\to} Co(g_*Q)$ which satisfy a 
transitivity condition. If now $Q$ is an equivcariant object of $\calC$, then
this family of isomorphisms makes $Co(Q)$ into an equivariant $\O^1_M$-torsor.

Then a curving $\na\in Co(P)\mapsto K(\na)$ is $G$-equivariant if and only
if and only if the relative $2$-forms $K(\na)$ and $K(\phi_*\na)$ coincide, for any
object $P$ of $[d_1^*\calC]_U$ and for any $\na\in Co(P)$.

\vskip 2.5pc

REFERENCES
\vskip .8pc

[{\bf A-M-M}] A. Alekseev, A. Malkin, E. Meinrenken, \it Lie Group Valued
Moment Maps, \rm  \ J. Differential Geom. \bf 48 \rm
(1998), no. 3, 445-495, dg-ga/9707021

[{\bf A-M-W}] A. Alekseev, E.Meinrenken, C. Woodward, \it Group-valued equivariant
localization, \rm dg/9905130

[{\bf B-D}] \jlb~ and  P. Deligne, \it Central extensions by $\un K_2$, \rm
preprint (1998), revised 2000

[{\bf B-M1}] \jlb ~and  D. McLaughlin, \it The geometry of the first
 Pontryagin class and of line bundles on loop
spaces\rm; I, 
 Duke Math. Jour. \bf 35 \rm (1994), 603-638 II,  Duke J. Math. \bf 83 \rm (1996), 105-139

[{\bf B-M2 }]  \jlb ~and  D. McLaughlin, \it The geometry of two-dimensional symbols, \rm
 K-Theory \bf 10 \rm (1996), 215-237

[{\bf Br1}] \jlb, \it Loop Spaces, Characteristic Classes and Geometric
Quantization, \rm Progress in Math. vol. \bf 107 \rm,
Birkha\"user (1993).

[{\bf Br2}] \jlb, \it Holomorphic gerbes and the Beilinson regulator,
  \rm  K-Theory (Strasbourg, 1992),
Ast\'erisque vol. 226 (1994), 145-174

[{\bf B-V}] N. Berline and M. Vergne, \it Z\'eros d'un champ de vecteurs et classes
caract\'eristiques \'equivariantes, \rm 
Duke Math. J. 50 (1983), no. 2, 539-549. 

[{\bf B-Z}] \jlb  ~and  B. Zhang, \it Equivariant K-theory of simply-connected Lie groups,
\rm  dg-ga/9710035

[{\bf Ch}] D. S. Chatterjee, \it On the construction of abelian gerbs, \rm D. Phil.
Cambridge Univ. (1998)

[{\bf D-D}] 
J. Dixmier and A. Douady, \it
Champs continus d'espaces hilbertiens et de $C^*$-alg\`ebres, \rm
Bull. Soc. Math. France \bf 91 \rm (1963), 227-284

[{\bf De}] P. Deligne,
\it Le symbole mod\'er\'e. \rm
 Publ. Math. IHES \bf 73, \rm (1991), 147-181. 

[{\bf D-F}] P. Deligne and D. Freed, \it
 Classical field theory. \rm in Quantum fields and strings: a course for mathematicians, Vol.
1 (Princeton, NJ, 1996/1997), 137--225, Amer. Math. Soc. (1999).\rm

[{\bf Gi}] J. Giraud, \it 
Cohomologie non ab\'elienne. \rm
 Grundlehren, Band 179. 
Springer-Verlag,  1971.

[{\bf G-H-J-W}]  
K. Guruprasad, J. Huebschmann, L. Jeffrey, and A. Weinstein, \it
Group systems, groupoids, and
moduli spaces of parabolic bundles, \rm  Duke Math. J. \bf 89 \rm (1997), no. 2, 377--412.

[{\bf Hi}] N. Hitchin, \it Lectures on Special Lagrangian Submanifolds, \rm DG/9907034

[{\bf K-S}] M. Kashiwara and P. Schapira,
\it Sheaves on manifolds, \rm corrected printing, 1994, \Sp

[{\bf P-S}]   A. Pressley and G. Segal, \it  Loop groups, \rm Oxford Univ. Press, New York,
1986.

[{\bf To}] V. Toledano Laredo, \it Positive energy representations of the loop groups
of non-simply connected Lie groups, \rm Comm. Math. Phys. \bf 207 \rm (1999) 2, 307-339

 [{\bf W-X}]  A. Weinstein and P. Xu, \it
 Extensions of symplectic groupoids and quantization, \rm 
J. Reine Angew. Math. 417 (1991), 159--189.

\vskip 2pc

Penn State University

Department of Mathematics

305 McAllister Bdg

University Park, PA 16802
\vskip .67pc

e-mail: jlb@math.psu.edu

\end